\pdfoutput=1
\documentclass[11pt]{article}
\usepackage{amssymb}
\usepackage{graphicx}
\parskip \baselineskip

\newcommand{\s}{\subseteq}

\newcommand{\HH}{\mathbb{H}}

\newcommand{\ol}{\overline}
\newcommand{\es}{\emptyset}

\newcommand{\Ra}{\Rightarrow}

\begin{document}

	\title{Compression with wildcards:\\
		All exact, or all  minimal hitting sets}
	
	\author{Marcel
		Wild\\[3pt]
		Department of Mathematical Sciences,
		University of Stellenbosch\\
		Private Bag X1, Matieland 7602, South
		Africa}
	
	\date{}
	\maketitle
	
	\begin{quote}
		A{\scriptsize BSTRACT}: {\footnotesize Our  objective is the {\it compressed} enumeration (based on wildcards) of all minimal hitting sets of general hypergraphs. To the author's best knowledge the only previous attempt towards compression, due to Toda [T], is based on BDD's and much different from our techniques. Traditional one-by-one enumeration schemes cannot compete  when the number of minimal hitting sets is large and the degree of compression is high. Our  method works particularly well in these two cases: Either compressing all {\it minimum cardinality} hitting sets, or compressing all {\it exact} hitting sets. }
	\end{quote}

	
	\section{Introduction}

	Let $W$ be a finite set (such as all sets in this article) and ${\cal P}(W)$ its powerset.  Given a hypergraph
	(=set-system) $\HH\subseteq {\cal P}(W)$, a ($\HH$-){\it hitting set} is a set $X\s W$ such that $X\cap H\neq \es$ for all {\it hyperedges} $H\in \HH$. Let ${\it HS}(\HH)$ be the set of all hitting sets, and ${\it MHS}(\HH)$ the subset of all (inclusion-)minimal hitting sets,
 henceforth called MHSes.
 The famous {\it Minimal Hitting Set Problem } is this: Given 
	$\HH\subseteq {\cal P}(W)$, is it possible to enumerate ${\it MHS}(\HH)$  in polynomial {\it total} time\footnote{An older synonym is {\it output polynomial time}.}, i.e. polynomial in  $w:=|W|,\ h:=|\HH|$, {\it and} $mhs:=|{\it MHS}(\HH)|$? We  refer to [1] and [2] for the history and the state of the art concerning this problem. 
	
	The  objective in our article is different and can be described in  picturesque ways as follows. For fixed $\HH$ identify the MHSes with diamonds and the ordinary hitting sets (i.e. the members of ${\it HS}(\HH)\setminus {\it MHS}(\HH)$) with worthless pebbles which, however, may be hard to distinguish from diamonds. Some friendly sponsor provides $R$ many nonempty boxes which  are filled with both kinds of stones.  {\it All} diamonds are distributed among the boxes but usually not all pebbles (which is just as well).  Our {\bf Main Quest} is to retrieve all diamonds (and only them) as efficiently as possible. A box is {\it good}
if it contains at least one diamond, and {\it bad} otherwise. A box 100\% filled with diamonds is {\it very-good}. As will be seen, depending on the structure of $\HH$,  very-good boxes can  both be numerous and heavy! Furthermore the number of diamonds in a very-good box is found at once, and the diamonds themselves are arranged in a pleasant, compressed  manner.

 To get a first impression of the quality of  boxes the Monte-Carlo method picks (say) 20 stones at random from each box $\rho$, and determines the number $\alpha(\rho)$ of diamonds among them. If $0<\alpha(\rho)<20$ then $\rho$ is {\it merely-good}, i.e. good but not very-good.  However, if $\alpha(\rho)=0$ then $\rho$ is only {\it likely-bad}, and if $\alpha(\rho)=20$ then $\rho$ is {\it likely-very-good}. If a likely-very-good row contains thousands of stones then classifying the stones one-by-one is time-consuming. Fortunately we will provide three  {\it criteria} for very-goodness which settle the issue faster.  Efficient criteria for badness are harder to come by but an elegant sufficient condition exists.
 As to merely-good boxes $\rho$, there are two approaches, each with benefits and drawbacks. The first is to classify the stones one-by-one. The second uses subtle machinery  but has the benefit that the diamonds in $\rho$ get {\it repackaged} into brandnew very-good boxes.
 
 Here comes the Section break-up, phrased in more mathematical terms. The preliminaries in Section 2 concern Boolean functions and three kinds of  wildcards; the $e$-, the $n$-, and the $g$-wildcard. All of them generalize  the  don't-care symbol $\ast$  familiar from describing partial models of Boolean functions. Furthermore we adopt the {\it Vertical Layout} technique used in data mining. In a nutshell, it substitutes set operations (e.g finding all suitable supersets of a given set) that involve {\it many small} sets by set operations with {\it few large}  sets. 
 Section 3 discloses the above-mentioned sponsor  (i.e. the transversal $e$-algorithm of [3]). Section 4 explains the mathematical nature of the $R$ boxes provided by the sponsor and  goes on (Theorem 1) to show that all minimum-{\it cardinality} MHSes occur in very-good boxes, which moreover can be pinpointed at once. Section 5 elaborates the first approach towards merely-good boxes by offering four  algorithms for one-by-one classification. Algorithm 1 relies on the diamonds (=MHSes) retrieved so far, whereas Algorithm 2 only relies on the knowledge of $\HH$. Algorithms 3 and 4 exploit  tricks that are fully justified only in Section 9. Section 6 elaborates the second approach towards merely-good boxes. Sections 7 and 8 propose two criteria (each of which sufficiend and necessary) for very-goodness. The first is based on inclusion-exclusion, the second on matroid theory (Rado's Theorem).
 
Section 9 introduces the key concept [4] of an {\it MC-set}. By definition  $X\s W$ is {\it MC} if for each $x\in X$ there is at least one hyperedge $H\in\HH$ that cuts $x$ sharply in the sense that 
$H\cap X=\{x\}$. The set-system ${\it MC}(\HH)$ of all MC-sets is dual to ${\it HS}(\HH)$ in that the former is a set-ideal, the latter a set-filter, and it holds (Theorem 2)  that ${\it MC}(\HH)\cap {\it HS}(\HH)={\it MHS}(\HH)$. Those subsets of $W$ which are not MC, yet all their proper subsets are  MC, are of particular importance. They are collected in the set-system ${\it MinNotMC}(\HH)$. For instance it allows us to calculate the cardinality $|{\it MHS}(\HH)|$ without knowing ${\it MHS}(\HH)$.
Section 10 calculates ${\it MinNotMC}(\HH)$. It exploits the fact that minimal set-coverings are cryptomorphic to minimal hitting sets and can hence be handled with the transversal $e$-algorithm. Section 11 features numerical experiments with Mathematica. 
 In a nutshell, our compression with wildcards works the better the fewer and the larger the hyperedges are.
In particular  very-good play a key role here. Although promising ideas  of previous Sections have not yet been implemented in Mathematica, in 11.6 we attempt a preliminary comparison of our methods with the algorithms of the two winners [4] and [5] of a competition carried out in [1].

Section 12 at first seems to abandon {\it minimal} hitting sets and turn to the  different topic of {\it exact} hitting sets (EHS). Is it that different? By definition $Y\s W$ is an EHS for $\HH$ if $|Y\cap H|=1$ for all $H\in \HH$. Under the mild assumption that $\bigcup\HH=W$ each EHS must be a MHS, yet the converse fails severly in that some hypergraphs have plenty MHSes and no EHSes. Nevertheless,  our previously  used $g$-wildcards can sometimes compress the set-system $E{\it HS}(\HH)$ of all hitting sets.  As to "sometimes", any fixed hypergraph $\HH\subseteq {\cal P}(W)$ induces a natural, apparently novel  equivalence relation $\sim$ on $W$. It turns out that compressing  $E{\it HS}(\HH)$ is possible iff $\sim$ is nontrivial. Furthermore
Knuth's popular Dancing-Link algorithm shows up in Section 12 and in Theorem 4 we enumerate the perfect matchings of any graph without $K_{3,3}$-minor in polynomial total time.

	\section{Preliminaries on  Boolean functions, partial models, wildcards, and Vertical Layout}

 After Boolean functions (2.1) we turn to $e$-wildcards (2.2-2.3), followed by $n$-wildcards and $g$-wildcards (2.4).  In 2.5 we  sieve the minimal members of any set-system ${\cal S}\s {\cal P}(W)$  and 2.6 introduces Vertical Layout.
 
  Throughout the article for any integer $w\ge 1$ we put  $ [w]:=\{1,2,..,w\}$. For convenience  usually  $W:=[w]$. If the powerset is concerned we write
 ${\cal P}[w]$ instead of ${\cal P}([w])$.
 Further we use the shorthand "iff' for "if and only if", and write $\subset$ (as opposed to $\s$) for proper inclusion.
	
{\bf 2.1}  We freely identify {\it bitstrings} of length $w$ (also called 01-{\it rows}) with subsets of
of $[w]$ in the usual way; thus $X=\{2,4,5\}$ (viewed, say, as subset of $[7]$) matches $x=(0,1,0,1,1,0,0)$. Depending on circumstances one or the other view is preferable.
 We now extend 01-rows to 012-{\it rows} such as \vspace{2mm}
 
 \centerline{$r=(0,2,2,1,0,2)$.}
 
 The following type of notation that refers to the positions of the various symbols will be used throughout:

 (1)\quad $zeros(r):=\{1,5\},\ ones(r):=\{4\},\ twos(r):=\{2,3,6\}.$

 While 01-rows encode  sets, 012-rows encode {\it set-systems} because '2' is viewed\footnote{In the literature often $\ast$ is used rather than 2.} as don't-care symbol which can be freely replaced by 0 or 1. Thus $r=(0,2,2,1,0,2)$ above encodes, and in fact will be {\it identified}\footnote{This is a bit sloppy but it outweighs the clumsiness of introducing an extra symbol for the represented set-system. From the context it will always be clear whether $r$ is meant to be a vector with entries from 0,1,2 or whether $r$ is a set-system.} with, the set-system
 
 $r=\{\{4\},\{4,2\},\{4,3\},\{4,6\},\{4,2,3\},\{4,2,6\},\{4,3,6\},\{4,2,3,6\}\}$

 which, with obvious shorthand notation (that will only be applied to sets of 1-digit numbers) can also be rendered (since elements of sets can be listed in arbitrary order) as

 $\{4,42,43,46,423,426,436,4236\}\ or\ as\ \{4,24,34,46, 234,246, 346, 2346\}.$

As to a general 012-row $r$, if it is viewed as a set-system, this  set-system is \\ $\{ones(r)\cup S:\ S\s twos(r)\}$. While $zeros(r)$ does not come up here, the 0's are as important as the 1's in the sequel (ponder what would become of $r=(0,2,2,1,0,2)$ without the 0's).

{\bf 2.1.1}  That leads us to $\{0,1\}$ viewed as Boolean algebra\footnote{We will only be concerned with the join and meet operations, so $1\vee 1=0\vee 1=1\vee 0=1,\ 0\vee 0=0$, and $0\wedge 0=0\wedge 1=1\wedge 0=0,\ 1\wedge 1=1$.} and to  {\it Boolean functions} $f:\{0,1\}^n\to\{0,1\}$ whose basic features are assumed to be familiar to the reader, so that we only need to fix notation here. Any bitstring $x\in\{0,1\}^n$ with $f(x)=1$ is a {\it model} of $f$. Apart from other means Boolean functions can be defined by {\it Boolean formulas}. Thus by writing $f(x):=x_1\vee x_2\vee x_3$ we define\footnote{In likewise fashion the formula defines a unique function $f:\{0,1\}^w\to\{0,1\}$ for every $w>4$. In the sequel it will always be clear which $w$ is meant. } a Boolean function $f:\{0,1\}^3\to\{0,1\}$ that e.g. satisfies
$f((0,1,1))=0\vee 1\vee 1=1$. It is clear that only $(0,0,0)$ fails to be a model, and so the {\it modelset} is

$Mod(f)=(2,2,2)\setminus\{(0,0,0,)\}=(1,2,2)\cup (2,1,2)\cup (2,2,1).$

The union on the righthand side is not disjoint since e.g. $(1,0,1)\in (1,2,2)\cap (2,2,1)$. Fortunately, this can be cured as follows (here and henceforth $\uplus$ signifies disjoint union):

\begin{center}
\begin{tabular}{cc}
		 $Mod(x_1\vee x_2\vee x_3)$ & \\
	= &	$({\bf 1},\ 2,\ 2)$ \\ 
        $\uplus$ & $(0,\ {\bf 1},\ 2)$ \\ 
        $\uplus$ & $(0,\ 0,\ {\bf 1})$  \\
	\end{tabular}
 \end{center}
		
This idea is long known and its visualization has been coined {\it Abraham-flag} in [6]. Thus a general $n\times n$ Abraham-flag has 1's in the main diagonal, 0's below it, and 2's above it. The row-cardinalities  sum up to $2^{n-1}+2^{n-2}+\cdots +1$ which equals $2^n-1$, as is to be expected. In connection with Boolean functions 012-rows usually describe {\it partial models}. For instance $(1,2,2)$ is a partial model of $x_1\vee x_2\vee x_3$ in the sense that replacing the 2's by 0 or 1 in any way results in a model of $x_1\vee x_2\vee x_3$.

{\bf 2.2} In addition to the don't-care symbol "2" we will use three further  wildcards. For starters, instead\footnote{However, Abraham-flags will reappear in 2.3 in new guise.} of using an $s\times s$ Abraham-flag to spell out $Mod(x_1\vee\cdots\vee x_s)$ we can, better still, simply define \vspace{2mm}

  \centerline{$ (e,e,...,e):=Mod(x_1\vee\cdots\vee x_s).$}

  Roughly speaking,  $s$ symbols $e$ (not necessarily adjacent) demand bitstrings to have "at least one 1 in that area". Combining such $e$-{\it wildcards} (distinguished by subscripts) gives rise to {\it 012e-rows} like

  (2)\quad  $r'=(e_1,0,2,e_1,e_2,1,0,e_2,2,2),$

which by definition consists of those subsets $S\s [10]$ that satisfy

\begin{itemize}
\item $2,7\not\in S$  (because $zeros(r')=\{2,7\}$)
\item $6\in S$ (because $ones(r')=\{6\}$)
\item $\{1,4\}\cap S=\emptyset$ (because of $e_1,e_1$)
\item $\{5,8\}\cap S=\emptyset$ (because of $e_2,e_2$)
\end{itemize}

The fact that $twos(r')=\{3,9,10\}$ reflects the fact that $3,9,10$ don't occur in any of the conditions. By {\it e-bubble} we mean the position-set of any given $e$-wildcard. Thus the $e_2$-bubble of the $e_2$-wildcard in (2) is $\{5,8\}$.
It is easy to see that \vspace{5pt}

\centerline{$|r'|=2^3\cdot (2^2-1)\cdot (2^2-1),$}

and that $2^2-1$ generalizes to $2^s-1$ for $e$-bubbles of size $s$. 

Alternatively (but clumsier) $r'$ in (2) could be defined\footnote{Provided $\{0,1\}^{10}$ is taken as domain of the Boolean function induced by the Boolean formula.} as

$(2')\quad r'=Mod(\ol{x_2}\wedge \ol{x_7}\wedge x_6\wedge (x_1\vee x_4)\wedge (x_5\vee x_8)).$

 {\bf 2.2.1} Observe that the intersection $\rho\cap\rho'$ of an 012e-row $\rho$ with an 012-row $\rho'$ is either empty (when 0's and 1's clash) or is again a 012e-row, which arises in obvious ways:

 \begin{center}
\begin{tabular}{cccccccccccc}
		& $\rho$ & $=$  & $(e_1,$ & $e_1,$ & $ e_1,$ & $e_2,$ & $e_2,$ &  $e_2,$ & $e_3,$ & $e_3,$ & $e_3)$ \\
      &  $\rho'$ & $=$  & $(2,$ & $2,$ & $0,$ &   $0,$ & $2,$ &  $2,$ & $1,$ & $2,$ & $0)$ \\
     & $\rho\cap\rho'$ & $=$  & $(e_1,$ & $e_1,$ & $ 0,$ & $0,$ & $e_2,$ &  $e_2,$ & $1,$ & $2,$ & $0)$ \\
	\end{tabular}
 \end{center}

{\bf 2.2.2} The set of all {\it minimal}\footnote{Recall that "minimal" always means minimal with respect to set inclusion.  } members contained in a 012e-row will play a crucial role. One checks that the set-system $Min(r')$ of all minimal members  of the set-system $r'$ in (2) equals

$(3)\quad Min(r')=\{615,618,645,648\}.$

Generally, if the 012e-row $r$ has $t\ge 1$ many $e$-wildcards of cardinalities $\epsilon_1,...,\epsilon_t$ then\footnote{For the special case of 012-rows $r$, i.e. $t=0$, we have $|Min(r)|=1$ and $deg(r)=|ones(r)|$.  } each $X\in Min(r)$ is of type $X=ones(r)\cup T$,  where $T$ cuts each $e$-bubble in exactly one element. Thus $|T|=t$. If we define the {\it degree} of $r$ as

$(4)\quad deg(r):=|ones(r)|+t,$

then

$(5)\quad Min(r)=\{X\in r:\ |X|=deg(r)\}\ and\ \quad |Min(r)|=\epsilon_1\cdot\epsilon_2\cdots\epsilon_t$.

For general set-systems ${\cal S}$ it will be more demanding (2.5) to sieve $Min({\cal S})$ from ${\cal S}$. Nevertheless (5) will keep coming back even in that context.

{\bf 2.3} Let us introduce higher-level Abraham-flags, i.e. constituted by certain 012e-rows as opposed to the 012-rows in 2.1. Consider

(6) $\quad r:=(e_1,e_1,e_2,e_2,e_3,\ e_1,e_2,e_2,e_3).$

Soon we need to be able  to e.g. sieve those bitstrings $(x_1,...,x_9)$ from $r$ that have at least one $1$ among $\{x_1,...,x_5\}$. In other words, we need to "impose" $(e,e,e,e,e)$ upon $r$, i.e. the {\it intersection} $r\cap (e,e,e,e,e,2,2,2,2)$ {\it of two 012e-rows} must be rewritten in a handy format. The answer is $r\cap (e,e,e,e,e,2,2,2,2)=r_1\uplus r_2\uplus r_3$ where

\begin{center}
\begin{tabular}{ccccccccccccc}
		(7)\quad & $r_1$ & $:=$  & $({\bf e_1},$ & ${\bf e_1},$ & $ e_2,$ & $e_2,$ & $e_3,$ & & $2,$ & $e_2,$ & $e_2,$ & $e_3)$ \\
      &  $r_2$ & $:=$  & $(0,$ & $0,$ & ${\bf e_2},$ &   ${\bf e_2},$ & $e_3,$ & & $1,$ & $2,$ & $2,$ & $e_3)$ \\
       &  $r_3$ & $:=$  & $(0,$ & $0,$ & $0,$ & $0,$ & ${\bf 1},$ & & $1,$ & $e_2,$ & $e_2,$ & $2)$ \\
	\end{tabular}
 \end{center}
 
The first part of the righthand side is a novel $3\times 3$ Abraham-flag in the sense that the boldface main diagonal entries are either $1$ (as in 2.1) or full $e$-wildcards. Likewise the entries below the main diagonal are again 0's. We leave it to the reader to figure out what happens above the main diagonal, and how all of this affects the last four columns in (7). See also Section 3.1.

{\bf 2.4} Dually to $e$-wildcards we will encounter  $n$-{\it wildcards} which demand "at least one 0 here". Thus for instance

  \begin{center}
\begin{tabular}{ccc}
		 $(n,n,n,n)$ &   $:=$   & $Mod(\ol{x_1}\vee \ol{x_2}\vee \ol{x_3}\vee \ol{x_4})$  \\
	               & =        &	$({\bf 0},\ 2,\ 2,\ 2)$ \\ 
     &   $\uplus$ & $(1,\ {\bf 0},\ 2,\ 2)$ \\ 
      &  $\uplus$ & $(1,\ 1,\ {\bf 0},\ 2)$  \\
     &   $\uplus$ & $(1,\ 1,\ 1,\ {\bf 0})$  \\
	\end{tabular}
 \end{center}

  Mutatis mutandis as in 2.2 we  define $n$-{\it-bubbles} and $012n$-{\it rows}.

 {\bf 2.4.1} Apart from $e$-wildcards and $n$-wildcards\footnote{We mention in passing that to some extent general clauses (i.e. with positive {\it and} negative literals) can  be handled  by mixing the two wildcards. For instance 
 $Mod(x_1\vee x_2\vee \ol{x_3}\vee \ol{x_4}\vee \ol{x_5})=(e,e,2,2,2)\uplus (0,0,n,n,n)$. Also in the present article the two wildcards will  appear simultaneously, if only in Section 9.},  a third type of wildcard takes care of the requirement "exactly one 1 here". Namely, by definition

 $(g,g,\ldots,g):=\{(1,0,...,0),(0,1,...,0),\ldots, (0,0,...,1)\}.$

 One trivial application of these $g$-{\it wildcards} (and coupled {\it g-bubbles}) is the compression of ${\it MHS}(\HH)$ for hypergraphs with disjoint hyperedges. Thus if $\HH_1=\{123,45,6789\}$ then  ${\it MHS}(\HH_1)=(g_1,g_1,g_1,g_2,g_2,g_3,g_3,g_3).$ Slightly more subtle and important later, one checks that $r'=(e_1,0,2,e_1,e_2,1,0,e_2,2,2)$ from (2) has $Min(r')=(g_1,0,0,g_1,g_2,1,0,g_2,0,0)$. Expressions like this are called {\it 01g-rows}.

	{\bf 2.5}  Let ${\cal S}\s {\cal P}([w])$ be any set system. The problem to get\footnote{All of the sequel applies mutatis mutandis to the set system $Max({\cal S})$ of all maximal members.} the set-system $Min({\cal S})$  of all minimal members of ${\cal S}$ occurs frequently in discrete mathematics. The naive way to proceed  is to decide for each $X\in {\cal S}$ whether there is another $Y\in {\cal S}$ with $Y\subset X$. Clearly $X$ belongs to $Min({\cal S})$ iff no such $Y$ exists. Since deciding whether or not $Y\subset X$ costs $O(w)$, the overall cost is $O(|{\cal S}|^2 w)$.
 
 To the author's best knowledge (readers are welcome to teach him better) the following  refinement has not appeared in the literature before. Start by grouping the members of ${\cal S}$ according to their cardinalities $m_1<m_2<\cdots <m_s$ (often  $m_{i+1}=m_i+1$). This induces the decomposition ${\cal S}={\cal S}[1]\uplus {\cal S}[2]\uplus\cdots\uplus {\cal S}[s]$. That costs $O(|{\cal S}|w)$. It suffices to show how to calculate $Min[i]:={\cal S}[i]\cap Min({\cal S})$ for all $1\le i\le s$. 
 
 Clearly $Min[1]={\cal S}[1]$ since minimum-cardinality implies minimal. Set  ${\cal S}'[i]:={\cal S}[i]$ for $2\le i\le s$.
 Throughout the remainder we will have $Min[i]\s {\cal S}'[i]\s {\cal S}[i]$ and the set-systems ${\cal S}'[i]$  keep shrinking until they reach ${\cal S}'[i]=Min[i]$. To begin with, pick any $X\in Min[1]$ and  remove all\footnote{This can be done "in one sweep" using the method of Vertical Layout discussed in 2.6.} $Y\in {\cal S}'[i]\ (i\ge 2)$ from ${\cal S}'[i]$ whenever $X\subset Y$. This costs $O(|{\cal S}|w)$.
 Doing the same for all  members $X'\in Min[1]$ costs $O(|{\cal S}|w\cdot |Min[1]|)=O(|{\cal S}|w\cdot min)$ where $min:=|Min({\cal S})|$. It is clear that afterwards ${\cal S}'[2]=Min[2]$. Next for each $X\in Min[2]$ and all $Y\in {\cal S}'[i]\ (i\ge 3)$ remove $Y$ from ${\cal S}'[i]$ whenever $X\subset Y$ (again VL can be used). Clearly afterwards ${\cal S}'[3]=Min[3]$. And so it goes on until we get ${\cal S}'[s]=Min[s]$.
 The overall cost is $O(|{\cal S}|w\cdot min\cdot s)=O(|{\cal S}|w^2\cdot min)$.

{\bf 2.6}  The operations $\vee,\ \wedge$ on $\{0,1\}$ extend to operations on $\{0,1\}^m$ (and they match union/intersection of sets in ${\cal P}([m])$). Adopting Mathematica terminology we call the extended operations $BitOr$ and $BitAnd$. For instance, referring to the columns of the $8\times 6$ matrix $A$  with rows $Z_1$ to $Z_8$ (Table 1), it holds that
	$BitAnd(col_2,col_6)=(0,0,1,1,0,0,0,1)^T$ (where the $T$ means 'transposed').

	\begin{tabular}{l|c|c|c|c|c|c|}
		& $col_1$ & $col_2$ & $col_3$ & $col_4$  &$col_5$  & $col_6$     \\ \hline
		& & & & & &   \\ \hline\hline
		$Z_1=$ &  1 & 1 & 1 & 0 & 0& 0  \\ \hline
		$Z_2=$ &  1 & 0 & 0 & 0 & 1& 0  \\ \hline
		$Z_3=$ &  1 & 1 & 0 & 0 & 0& 1  \\ \hline
		$Z_4=$ &  0 & 1 & 0 & 0 & 1& 1  \\ \hline
		$Z_5=$ &  1 & 0 & 1 & 1 & 0& 0  \\ \hline
		$Z_6=$ &  0 & 0 & 1 & 1 & 1& 0  \\ \hline
		$Z_7=$ &  0 & 0& 1 & 1 & 0& 1  \\ \hline
		$Z_8=$ &  0 & 1 & 0 & 1 & 0& 1  \\ \hline\hline
	\end{tabular}
	
	{\sl Table 1: Illustrating Vertical Layout.}
	
{\bf 2.6.1}	What is this good for?  The fact that $BitAnd(col_2,col_6)$  had a component 1 exactly on the 3th, 4th and 8th position tells us that among the sets $Z_1,...,Z_8$ the ones that contain the set $\{2,6\}$ are exactly $Z_3,Z_4,Z_8$.  This is e.g. relevant for speeding up the method of 2.5.
 
 {\bf 2.6.2} Here comes another application. Consider  the set system
	
	$(8)\quad {\cal G}:=\{\{1,2,3\},\{1,5\},\{1,2,6\},\{2,5,6\},\{1,3,4\},\{3,4,5\},\{3,4,6\},\{2,4,6\}\}.$
	
	The straightforward (='horizontal') way to see whether $X=\{1,2,5\}$ is a ${\cal G}$-transversal checks whether any intersection $X\cap Y\ (Y\in {\cal G})$ is empty. In contrast, {\it Vertical Layout (VL)} demands\footnote{For the history of VL see e.g. arXiv:2002.09707.}  to build the $8\times 6$ matrix $A({\cal G})$ whose $i$th row $Y_i'$ is the characteristic bitstring of the $i$th set $Y_i$ listed in (8). It happens that $A({\cal G})$ is rendered in Table 1. A moment's reflection confirms the following. The fact that
	$BitOr(col_1,col_2,col_5)=(1,1,1,1,1,1,0,1)^T$
	does not equal $ (1,1,1,1,1,1,1,1)^T$, is tantamount to $X$ not being a ${\cal G}$-hitting set ($X\cap Y_7=\es$).  Although the formal complexities of the horizontal and vertical way  coincide, in practise VL is the faster the more (small) sets ${\cal G}$ contains. Simply put, computer hardware prefers doing few operations with long bitstrings to doing many operations with short bitstrings.

\section{Review  of the transversal $e$-algorithm}

We survey the transversal $e$-algorithm (3.1) and adapt it to count or generate  hitting sets of fixed cardinality (3.2). In 3.3 we indicate how the transversal $e$-algorithm dualizes to the noncover $n$-algorithm.

{\bf 3.1} Consider the task to enumerate the set ${\it HS}(\HH_2)$ of all hitting sets of the hypergraph $\HH_2$ whose five hyperedges $X\s [6]$ are

 $(9)\ H_1=\{1,2,5\},\ H_2=\{3,4\},\ H_3=\{4,5,6\},\ H_4=\{1,3,5\},\ H_5=\{2,6\}.$

 One idea is to first compute the hitting sets of the hypergraph $\{H_1\}$, then the ones of $\{H_1,H_2\}$, and so forth until we obtain the hitting sets of $\{H_1,...,H_5\}=\HH_2$. Calculating ${\it HS}(\{H_1\})$ is easy in view of 2.2. It consists of all bitstrings (=subsets of [6]) that have at least 1 on the positions 1,2,5, and so ${\it HS}(\{H_1\})=(e,e,2,2,e,2)$. Likewise
 ${\it HS}(\{H_1,H_2\})=(e_1,e_1,e_2,e_2,e_1,2)=:r'$.
 
 Now it gets trickier because $H_3$ intersects $H_1$ and $H_2$, i.e. the $e_3$-wildcard supposed to be modeling $H_3$ interferes with existing $e$-wildcards. In 2.3 we indicated how this is to be handled. Recall that the row in (6), which suffered the same predicament as $r'$ above, had to be replaced by three {\it candidate sons} in (7). The essence of the transversal $e$-algorithm is to keep on picking the topmost row  $r'$ of a "to do" 
  stack of 012e-rows and to impose some $e$-wildcard upon $r'$, which in turn can trigger up\footnote{Here $t$ is as in (4) and (5).
  Concerning the "to do" stack, the standard name is Last-In-First-Out (LIFO) stack.
   LIFO-stacks are standard data structures which match the depth-first search of trees.} to $t$ candidate sons.
 Each candidate son $r_i$ must be {\it feasible} in the sense that $r_i\cap {\it HS}(\HH)\neq\es$, for otherwise further processing of $r_i$ cannot possibly yield any hitting sets. The feasible candidate sons are put on top of the LIFO stack, the others are discarded. Fortunately deciding feasibility is easy:

 (10) $r$ {\it is feasible iff} $(\forall H\in\HH)(H\not\s zeros(r)).$

The effect of discarding infeasible candidate sons is that in each set of candidate sons at least one will be feasible. This in turn is the reason that the $e$-algorithm runs in total polynomial time, in fact in $O(Rh^2w^2)$ time. For the fine details of this {\it transversal} $e$-{\it algorithm}\footnote{Due to its use in previous publications we stick with 'transversal e-algorithm'. Other than that we always use "hitting set" instead of the synonym "transversal".}   the reader is referred to [3].
 To summarize, for any given hypergraph $\HH\s {\cal P}([w])$  the transversal $e$-algorithm  renders ${\it HS}(\HH)$ as a disjoint union of $R$ many 012e-rows, thus

 (11)\quad ${\it HS}(\HH)=\biguplus_{i=1}^R \ol{\rho_i}.$

 {\bf 3.1.1} Applied to $\HH_2$ the transversal $e$-algorithm yields ${\it HS}(\HH_2)=\ol{\rho_1}\uplus\cdots\uplus \ol{\rho_4}$, where the $\ol{\rho_i}$'s are defined in Table 2.

 \begin{tabular}{l|c|c|c|c|c|c|}
		\hline
		$\ol{\rho_1}=$ &  $e$ & $e$ & 1 & 0 & 0& 1 \\ \hline
		$\ol{\rho_2}=$ &  2 & $e_1$ & $e_2$ & $e_2$ & 1& $e_1$  \\ \hline
		$\ol{\rho_3}=$ &  0 & 1 & 1 & 1 & 0& 2  \\ \hline
		$\ol{\rho}_4=$ &  $1$ & $e$ & 2 & 1 & 0& $e$ \\ \hline	
	\end{tabular}
	
	{\sl Table 2: Representing ${\it HS}(\HH_2)$ as disjoint union of  $012e$-rows}

 In view of 2.2 we conclude that

 $|{\it HS}(\HH_2)|=|\ol{\rho_1}|+\cdots |\ol{\rho_4}|= (2^2-1)+2(2^2-1)^2+2+2(2^2-1)=29.$

{\bf 3.2} Let $\mu:=\mu(\HH)$ be the minimum cardinality achieved by any hitting set of the hypergraph $\HH$.  Often $\mu$ gets known\footnote{According to [8] the cost of finding a minimum-cardinality transversal is $O(1.2381^n)$ where $n$ is the sum of $w$ and all cardinalities $|H|\ (H\in \HH)$.} only {\it after} (11) has been obtained. For all $c\in \{\mu,\mu+1,..,w\}$ we put

 \begin{itemize}
 \item[(12)] ${\it HS}(\HH,c):=\{X\in {\it HS}(\HH):\ |X|=c\}.$
 \end{itemize}

 Of particular interest is the set-system 

 \begin{itemize}
 \item[(13)] ${\it MCHS}(\HH):={\it HS}(\HH,\mu)\s MSH(\HH).$
 \end{itemize}
	
{\bf 3.2.1} In some circumstances (e.g. in [7]) it is irrelevant whether the $\HH$-hitting sets are minimal; just their cardinality matters. Let us hence calculate
 $|{\it HS}(\HH,c)|$ for any fixed $c\ge\mu$. Viewing (11) for any such $c$ let $\ol{I}$ be the set of indices $i\le R$ such that the 012e-row $\ol{\rho_i}$ has degree $\le c$. (That's because  $\ol{\rho_i}\cap {\it HS}(\HH,c)=\es$ if $i\not\in\ol{I}$.) Putting $S(\ol{\rho_i}):=\{X\in \ol{\rho_i}:\ |X|=c\}$ we get $|{\it HS}(\HH,c)|$ by summing up the numbers
$|S(\ol{\rho_i})|\ (i\in\ol{I})$. It is easy to calculate the numbers $|S(\ol{\rho_i})|$ with inclusion-exclusion; for a faster way see [3,Thm.1].

{\bf 3.2.2} Suppose the set ${\it HS}(\HH,c)$ {\it itself} needs to be calculated. By the above each fixed set-family ${\it HS}(\HH,c)$ is the disjoint union of all sets $S(\ol{\rho_i})\ (i\in\ol{I})$. But  sieving $S(\ol{\rho_i})$  from  $\ol{\rho_i}$ is more cumbersome than  calculating
$|S(\ol{\rho_i})| $. Leaving ways of compression to the future, we only note that	if $S(\ol{\rho_i})$ has $\alpha$ elements then by [3,Thm.2] it can be enumerated one-by-one in total polynomial time $O(\alpha w^2)$.

{\bf 3.2.3} If  ${\it HS}(\HH,c)$ is of interest cardinality-wise (or the members themselves) for {\it all} values $\mu\le c\le w$, then upon  running the transversal $e$-algorithm each $c$ gets processed as discussed in 3.2.1 (or 3.2.2).
However, if only values  $c\le d$ for some {\it bound} $d$ are relevant, then it pays to {\it adjust} the transversal $e$-algorithm as follows. In addition to (10), the arising candidate sons should also satisfy $deg(\ol{\rho})\le d$. That's because $deg(\ol{\rho})>d$ implies that all successor rows $\ol{\rho_i}$ of $\ol{\rho}$ will have $deg(\ol{\rho_i})\ge deg(\ol{\rho})>d$, and so cannot contain any members of ${\it HS}(\HH,c)$. 
Problem is, in contrast to the remarks after (10) it can now happen that  some rows loose {\it all} their candidate sons. Nevertheless, performance in practise may be good.

 {\bf 3.3} The family ${\it HS}(\HH)$ of all $\HH$-hitting sets  is a {\it set-filter} ${\cal F}$ in the sense that\\  $(X\in {\cal F}\ and\ X\s Y)\Ra Y\in {\cal F}$. Now let ${\cal S}\s {\cal P}([w])$ be a set system. Call $Z\in {\cal P}([w])$ a ${\cal S}$-{\it noncover} if $Z\not\supseteq Y$ for all $Y\in {\cal S}$. Then the family $NC({\cal S})$ of all ${\cal S}$-noncovers is a {\it set-ideal} ${\cal J}$ in the sense that $(X\in {\cal J}\ and\ Y\s X)\Ra Y\in {\cal J}$. 
Consider {\it any} set-filter ${\cal F}\s {\cal P}([w])$. The minimal members of ${\cal F}$ are called its {\it generators} and they determine ${\cal F}$ uniquely.  Likewise for {\it any} set-ideal\footnote{Set-ideals are also called (abstract) {\it simplicial complexes.}} ${\cal J}\s {\cal P}([w])$ the maximal members of ${\cal J}$ are called its {\it facets} and they determine ${\cal J}$ uniquely. Furthermore, let ${\cal F}$ and ${\cal J}$ be complementary set-systems in the sense that
 ${\cal F}\uplus{\cal J}={\cal P}([w])$. It then holds that ${\cal F}$ is a set-filter iff ${\cal J}$ is a set-ideal.
 

  Given $\HH\s {\cal P}([w])$ the transversal $e$-algorithm renders the set-filter ${\it HS}(\HH)$ in the convenient format (11). Since set-filter and set-ideal are dual concepts, and so are $e$-wildcards and $n$-wildcards, it comes as no surprise  that some {\it noncover n-algorithm} (see e.g. [6]), fed with ${\cal S}$
 renders the set-ideal $NC({\cal S})$ as a disjoint union of $R'$ many 012n-rows:

 (11')\quad $NC({\cal S})=\biguplus_{i=1}^{R'} \ol{\sigma_i}.$

	\section{From minimum-cardinality toward inclusion-minimal}

 We show that (11) persists even when all 012e-rows $\ol{\rho_i}$ get "shaved" and become certain 01g-rows $\rho_i\s\ol{\rho_i}$.
 Thus (11) improves to (17). It turns out  that $MC{\it HS}(\HH)$ in (13) is the union of some such rows ${\rho_i}$. 
 In 4.2 we comment on situations where $MC{\it HS}(\HH)={\it MHS}(\HH)$, and in 4.4 resume the Monte Carlo of Section 1 in order to get an estimate for $|{\it MHS}(\HH)|$..

 {\bf 4.1} Let $\HH\s {\cal P}([w])$ be a hypergraph. In the remainder of the article we assume that the transversal $e$-algorithm has rendered ${\it HS}(\HH)$ as a disjoint union of $R$ many 012e-rows as in (11). Different from [3] where these rows were coined 'final', here the availability of them is not the end but only the beginning.  That's why we henceforth call them \vspace{2mm}
 
\centerline{\it semifinal 012e-rows.} 
 
 Suppose $X\s [w]$ is any minimal $\HH$-hitting set. Then $X$ is contained in some  semifinal 012e-row $\ol{\rho_i}$ because of (11). Being minimal within ${\it HS}(\HH)$, a fortiori $X$ is minimal within the smaller set-system
  $\ol{\rho_i}\s {\it HS}(\HH)$, i.e.  $X\in Min(\ol{\rho_i})$. In view of (5) it follows that for all $1\le i\le R$:

  $(14)\quad  \ol{\rho_i}\cap {\it MHS}(\HH)\s Min(\ol{\rho_i})=\{X\in \ol{\rho_i}:\ |X|=deg(\ol{\rho_i})\}.$

  In particular consider $Y\in MC{\it HS}(\HH) \s {\it MHS}(\HH)$. As before $Y\in Min(\ol{\rho_j})$ for some $j\le R$. But all sets in $Min(\ol{\rho_j})$ have the same cardinality as $Y$, and so are themselves in $MC{\it HS}(\HH)$. Hence $\s$ in (14) becomes =.
   To summarize:

  {\bf Theorem 1: }{\it Assume that ${\it HS}(\HH)$ is represented as disjoint union of 012e-rows $\ol{\rho_i}$ as in (11). Then, with notation as above, $MC{\it HS}(\HH)$ is the disjoint union of those sets $Min(\ol{\rho_i})$ that have $deg(\ol{\rho_i})=\mu$.}

  To illustrate consider ${\it HS}(\HH_2)=\ol{\rho_1}\uplus\cdots\uplus\ol{\rho_4}$ in Table 2.
  One checks that all these rows happen to have degree $3$, and so $\mu=3$. It follows from Theorem 1 and the fact (see 2.4.1) that sets of type $Min(\ol{\rho_i})$ can conveniently be rendered by single $01g$-rows $\rho_i$ that 

  $(15)\quad MC{\it HS}(\HH_2)=\rho_1\cup \rho_2\cup \rho_3\cup \rho_4,$
 
 where the $\rho_i$'s are defined below:

 \begin{itemize}
 \item[(16)] $Min(\ol{\rho_1}))=\{136,236\}=(g,g,1,0,0,1)=:\rho_1$
  \item[] $Min(\ol{\rho_2}))=\{235,245,356,456\}=(0,g_1,g_2,g_2,1,g_1)=:\rho_2$
   \item[] $Min(\ol{\rho_3}))=\{234\}=(0,1,1,1,0,0)=:\rho_3$
    \item[] $Min(\ol{\rho_4}))=\{124,146\}=(1,g,0,1,0,g)=:\rho_4$
 \end{itemize}

{\bf 4.2} As opposed to (15) where incidently $MC{\it HS}(\HH_2)={\it MHS}(\HH_2)$, for general hypergraphs $\HH$ only few semifinal 012e-rows $\ol{\rho_i}$ will have degree $\mu$!  If only $MC{\it HS}(\HH)$ is sought then all rows $\ol{\rho_i}$ with $deg(\ol{\rho_i})>\mu$ are superfluous. Yet to avoid them one cannot  proceed as in 3.2.3 because usually $d:=\mu$ is not known in advance.
However, guessing and working with some slightly larger $d>\mu$ will still beat computing all $R$ rows. (If it happens that one guesses a $d<\mu$ then the proposed method will not deliver any semifinal 012e-rows. But it  will improve the next guess, and with binary search one can even pin down $\mu$.)

{\bf 4.2.1} Interestingly, in the following set-up $\mu$ {\it is} known\footnote{Readers aware of  other scenarios of that type, please  let the author know.} in advance; it even happens  that ${\it MHS}(\HH)=MC{\it HS}(\HH)$. Namely, if $\HH$ is the family of all cocircuits [9,p.653] of a matroid then ${\it MHS}(\HH)$ is the set of all matroid bases and $\mu$ is easy to come by. In arXiv:2002.09707 (submitted) this has been implemented for the scenario where the cocircuits are the minimal cutsets of a graph $G$, in which case ${\it MHS}(\HH)=MC{\it HS}(\HH)$ is the set of all spanning trees of $G$. 

{\bf 4.2.2} Suppose that $\mu$ is known, be it by binary search or by theoretical reasoning as in 4.2.1. Then one still sits with the problem (mentioned in 3.2.3) that some top-rows of the LIFO stack may loose all their candidate sons. That this {\it cannot} happen in 4.2.1 is one of the (numerically well-supported) conjectures raised in arXiv:2002.09707. In another vein, if all $H\in\HH$ have $|H|=2$, so $\HH$ is the edge-set of a graph, then instead of MHSes one rather speaks of minimal vertex-covers. In this scenario $\mu$ remains hard to compute, but at least "loosing all candidate sons" can be avoided (work in progress).

{\bf 4.3} Generalizing Table 2 and (16), each  semifinal 012e-row $\ol{\rho}:=\ol{\rho_i}$ appearing in (11)  yields the 

$\centerline{\it semifinal 01g-row (or simply: semifinal row)} $ 

 $\rho:=Min(\ol{\rho})$ where all 2's of $\ol{\rho}$ have been replaced by 0's and each $e$-wildcard of length $\epsilon_j$ has been replaced by  a $g$-wildcard of the same length $\gamma_j:=\epsilon_j$. Hence, akin to (4) and (5), the semifinal 01g-row $\rho$ has $t$ many $g$-wildcards and it holds that the $\gamma_1\gamma_2\cdots\gamma_t$ members of $\rho$ all have cardinality $|ones(\rho)|+t$.
It follows from (11) that

$(17)\quad {\it MHS}(\HH)\s \biguplus_{i=1}^R \rho_i =:SF(\HH).$

Accordingly we have

$(17')\quad mhs=|{\it MHS}(\HH)|\ \le\ sf:=| \biguplus_{i=1}^R \rho_i|.$

The following terminology will be handy as well.
A   semifinal 01g-row $\rho$ is {\it bad} if ${\it MHS}(\HH)\cap\rho=\emptyset$, and {\it good} otherwise. Additionally call $\rho$  {\it very-good} if $\rho\s {\it MHS}(\HH)$, and call $\rho$  {\it merely-good} if it is good but not very-good. Each $X\in \rho\setminus {\it MHS}(\HH)$  a {\it dud}. For 012-rows it holds that $good\Leftrightarrow$ $very$-$good$.

{\bf 4.4} A simple attempt to settle  "good or bad?" is the {\it Monte-Carlo} way. That is, pick uniformly and at random  $X\in\rho$ and check (in whatever way) whether or not $X\in MHS$. If yes, then $\rho$ is good. If no,
test some more $X$. The more often the answer persists to be no, the likelier  $\rho$ is bad. As to "likelier", the density $d:=|\rho\cap {\it MHS}(\HH)|\ /\ |\rho|$ can be estimated to any desired precision $\epsilon$ as follows. Given $\epsilon,\delta>0$, standard statistics yields  a value $d'$ such that 
(with error-probability $<\delta$) it holds that $d\in [(1-\epsilon)d',(1+\epsilon)d']$. Since $|\rho |=\gamma_1\cdots\gamma_t$ is known,  $d'$ also yields an estimate for $|\rho\cap {\it MHS}(\HH)|$, and hence in view of (17) for $|{\it MHS}(\HH)|$ .

 \section{Four ways to sieve the MHSes from the semifinal 01g-rows }

 Let $\rho$ be a fixed semifinal 01g-row. In this Section we present four methods (Algorithm 1 to Algorithm 4) to {\it classify} all $X\in\rho$ one-by-one, i.e. to decide whether $X$ is an MHS or a dud. Algorithm 1 relies on 2.5 and 2.6.1, whereas Algorithm 2 uses the kind of Vertical Layout  in 2.6.2. Algorithms 3 and 4 rely on presently  "magic" set-systems $PotKi(\rho)$, respectively $ {\it MC}(\HH)$, whose capabilities and whose acquiry will be postponed to later Sections. 

 {\bf 5.1} Referring to 2.5 let $m_1<m_2<\cdots<m_s$ be the numbers that occur as cardinalities of $\HH$-hitting sets. Then\footnote{Often, yet not always, $m_{i+1}=m_i+1$ for all $1\le i\le s-1$.  } $m_1=\mu$ and $m_s=w$. Putting ${\cal S}:=SF(\HH)$ we have $Min({\cal S})={\it MHS}(\HH)$, and following 2.5 we get $MSH(\HH)$ in time $O(sf\cdot w\cdot mhs\cdot s)$. This
{\it Algorithm 1} may compare favorably to methods in spe if $s$ gets replaced (say) by 3, i.e. if we only care for MHSes of cardinality $\le \mu+2$. In this case the cost is $O(sf\cdot w\cdot mhs)$.

{\bf 5.2}  Let us view the hyperedges of $\HH$ as bitstrings and take them as the rows of an $h\times w$ matrix $A$. Fix a semifinal 01g-row $\rho$ and put $k:=deg(\rho)$. For each fixed $Y\in \rho$ (hence a hitting set) it holds that $Y\in {\it MHS}(\HH)$ iff {\it no} set $X:=Y\setminus\{a\}\ (a\in Y)$ is an hitting set. Whether or not VL based on $A$ (see 2.6) is used, the formal cost to classify $X$ is $O(kh)$. Hence classifying $Y$ costs $O(k^2w)$.
Furthermore finding $\rho\cap {\it MHS}(\HH)$ costs $O(|\rho|k^2h)$, and finding ${\it MHS}(\HH)$ with the sketched  {\it Algorithm 2} costs $O(sf\cdot w^2h)$ (since $|Y|=k$ becomes $|Y|\le w$).

{\bf 5.2.1} Observe that the bound  in 5.1 to get ${\it MHS}(\HH)$ does not adapt smoothly to a bound for calculating
just $\rho\cap {\it MHS}(\HH)$. This contrasts with 5.2 where we obtained $O(|\rho|k^2h)$ for the latter task, due to the fact that it costs $O(k^2h)$ to decide whether any $k$-hitting set  of a hypergraph with $h$ hyperedges is minimal. While the $O(k^2h)$ bound for this basic decision problem has probably been observed before, let us indicate a  surprising improvement of it. To fix ideas, suppose that $k=5$ and that the minimality of an $\HH$-hitting $Y=\{1,2,3,4,5\}$ (where $|\HH|=h$) needs to be decided. In 2.6.2 the VL way to handle $Y$ demands to calculate

$col_{1234}:=col_1\vee col_2\vee col_3\vee col_4\ and\ col_{1235}:=col_1\vee col_2\vee col_3\vee col_5$ 

$and\ col_{1245}:=col_1\vee col_2\vee col_4\vee col_5\ and\ col_{1345}:=col_1\vee col_3\vee col_4\vee col_5$

$and\ col_{2345}:=col_2\vee col_3\vee col_4\vee col_5.$

This required $5\cdot 3=15$ basic BitOr operations, but one can improve that to 11:

$col_{12}:=col_1\vee col_2,\ col_{123}:=col_{12}\vee col_3,\ col_{124}:=col_{12}\vee col_4,\ col_{1234}:=col_{123}\vee col_4,$

$col_{1235}:=col_{123}\vee col_5,\ col_{1245}:=col_{124}\vee col_5,\ col_{45}:=col_4\vee col_5,\ col_{345}:=col_3\vee col_{45},$

$col_{245}:=col_2\vee col_{45},\
col_{1345}:=col_1\vee col_{345},\ col_{2345}:=col_3\vee col_{245}.$

 Driving this idea further\footnote{Interested readers are welcome to help refine the author's handwritten notes into publishable form.} one can improve $O(k^2w)$ to $O(k^{4/3}w)$.

 {\bf 5.3} Given a semifinal 01g-row $\rho$ suppose it was possible (more on that in 9.3) to get a set-system $PotKi(\rho)\s {\cal P}([w])$  such that any given $X\in\rho$ is a dud iff it gets {\it killed} by some $Z\in PotKi(\rho)$ in the sense that $Z\s X$. So suppose the toy row $\tilde{\rho}:=(1,g_1,g_1,g_1,g_2,g_2)$ has $PotKi(\tilde{\rho})=\{15,126\}$. Since $15$ kills $152,153,154$ and $126$ kills $126$, we have four duds and hence $\tilde{\rho}\cap {\it MHS}(\HH)=\{136,146\}$.
The availability of $PotKi(\rho)$ facilitates a lot  the calculation of 

\centerline{$Duds(\rho):=\{X\in\rho:\ X\ is\ a\ dud\}\ (=\rho\setminus {\it MHS}(\HH))$.}

Namely, embarking onto VL (which makes the more sense the larger $\rho$), we view the members of $\rho$ as bitstrings and take them as the rows of a $|\rho|\times w$ matrix $A$. Starting with $Duds(\rho):=\es$ we process $PotKi(\rho)$ one by one and update $Duds(\rho)$ accordingly as follows {\it (Algorithm 3)}. Say $Z=\{2,4,7\}\in PotKi(\rho)$. If $col_i$ is the $i$-th column of $A$ we calculate $col:=col_2\wedge col_4 \wedge col_7$. Then $ones(col)$ is the set of row-numbers whose corresponding rows of $A$ get killed by $Y$. Thus we update $Duds(\rho)=Duds(\rho)\cup ones(col)$.

{\bf 5.3.1} Picked from the author's random experiments, here comes a more demanding semifinal $01g$-row $\rho$. It is defined by
$zeros(\rho):=\es,\ ones(\rho):=\{4,5,6\}$,  has $g$-bubbles\\  $\{1,8\},\{2,10\},\{3,11\},\{7,9,12\}$, and has 

$PotKi(\rho)=\{Z_1,Z_2,Z_3,Z_4\}:=\{\{2,4,6\},\{2,5,6\},\{1,5,10\},\{5,8,10\}  \}$.
	
	Here $Z_1$ kills (exactly) the twelve  sets of type $\{{\bf 4},5,{\bf 6}\}\cup\{a,{\bf 2},c,d\}$, $Z_2$ kills the  sets of type $\{4,{\bf 5},{\bf 6}\}\cup\{a,{\bf 2},c,d\}$, i.e. the
	{\it same} as before, $Z_3$ kills the six sets of type $\{14,{\bf 5},6\}\cup\{{\bf 1,10},c,d\}$, and $Z_4$ the six sets of type $\{4,{\bf 5},6\}\cup\{{\bf 8, 10},c,d\}$.       
	Since the killed sets happen to be either identical or disjoint, it follows from $12+6+6=|\rho|$ that  $\rho$ gets killed entirely.
 It is an example of a 'sophisticated-bad' row, the exact definition following in 9.1.

{\bf 5.4} Fix some hypergraph $\HH\subseteq {\cal P}([w])$. Following  [4] we say  $S\subseteq [w]$ is a {\it MC-set} (or: {\it is MC}) iff for each $b\in S$ at least one $H\in\HH$ {\it cuts $b$ sharply from } $S$, i.e. $H\cap S=\{b\}$. Each other $H'\in\HH$ either cuts out $b$ sharply as well, or has $|H'\cap S|\ge 2$,  or has $H'\cap S=\es$.  It is evident that a subset of a MC-set is again a MC-set. Hence the family 

(18)\quad ${\it MC}(\HH):=\{S\subseteq [w]:\ S\ is\ MC\}$ {\it is a set-ideal.}

 Note that MC-sets need not be hitting sets. To witness take $\HH_3:=\{\{1,3\},\{2,4\},\{3,4\}\}$.
	One checks that $\{1,2\}$ is a  MC-set yet not an $\HH$-hitting set. However, it holds (Section 9) that  $\rho\cap {\it MC}(\HH)=\rho\cap {\it MHS}(\HH)$ for each semifinal 01g-row $\rho$. 
 
 This suggests an elegant method for classifying any $X\in \rho$.
Namely, initialize a {\it testset} $T$ to $T:=\es$. Process all $H\in \HH$ and update $T:=T\cup (X\cap H)$ (programmer's speak) whenever $|X\cap H|=1$. As soon as $T=X$ occurs, we know that $X\in {\it MHS}(\HH)$. If it never occurs then $X\not\in {\it MHS}(\HH)$. Since classifying $X$ that way costs
$O(hw)$ we have a method, call it {\it Algorithm 4}, that calculates ${\it MHS}(\HH)$ in time $O(sf\cdot hw)$.

{\bf 5.4.1} Let us indicate how VL may further speed up calculating $\rho\cap {\it MHS}(\HH)$ (as always, without challenging the formal bound, in this case $O(|\rho|hw)$). For starters, the $wh$ many sets

$S(a,H):=\{Z\in\rho:\ a\in H\cap Z,\ |H\cap Z|\ge 2\}\hspace{1cm}  (a\in [w],\ H\in\HH)$

need to be calculated. To do so initialize all of them as $S(a,H):=\es$. Next for each fixed $Z\in\rho$ do the following. Using VL determine all $K\in\HH$ with $|K\cap Z|\ge 2$. For any such $K$ add\footnote{Instead of adding whole sets $Z$, code the sets as {\it numbers} $f(Z)$. Because $\rho$ is a 01g-row, there is a natural "lexicographic" encoding. To witness, take $\rho=(g_1g_1,g_2g_2,g_3g_3)$ which for simplicity has $zeros(\rho)=ones(\rho)=\es$. Put $f(\{2,4,6\})=f(246):=1,\ f(245):=2,\ f(236):=3,...,f(135):=8$.} $Z$ to all sets $S(a,K)\ (a\in K\cap Z)$. For $a\in [w]$ and $Z\in\rho$ with $a\in Z$ call 
$Z$ an {\it $a$-dud} if there is no hyperedge that sharply cuts out $a$ from $Z$ (and so $Z\not\in {\it MHS}(\HH)$). It is easy to see that

$S(a):=\bigcap\{S(a,H):\ a\in H\in\HH\}$

is the set of all $a$-duds, and that VL speeds up the calculation of $S(a)$ the more the bigger $|\rho|$. Consequently 

$Duds:=\bigcup\{S(a):\ a\in [w]\}$ 

is the set of all duds contained in $\rho$. Put another way, $\rho\cap {\it MHS}(\HH)=\rho\setminus Duds$.

	\section{Replacing merely-good rows by  very-good rows}
	
	In Section 5 we presented four algorithms to unravel the MHSes contained  in a fixed semifinal 01g-row $\rho$. Any such MHS, viewed as bitstring $x\in \{0,1\}^w$ is\footnote{To be pedantic, since by definition every 012g-row is a {\it set} of bitstrings, not $x$ but $\{x\}$ is a very-good row.} is a very-good row, and so one could say that each semifinal row $\rho$ is either bad or can be represented as a disjoint union of very-good rows. But it would  be nice to use {\it fewer} than $|\rho|$ very-good rows to exhaust $\rho$. 
 
 Suppose we possess (more on that later) criteria that allow us to quickly classify each semifinal 01g-row as bad, merely-good, and very-good.  The bad ones are thrown away, the very-good ones are in optimal shape, but what about the merely-good rows $\rho$? Aren't we back to Square 1 and need to scan $\rho$ one by one?. Not so.  We start with a toy example in 6.1 and follow up with theory in 6.2.

	{\bf 6.1} Consider the hypergraph $\HH_4\s{\cal P}([6])$ with hyperedges
	
	(19)\quad $H_1=\{1,5,6\},\ H_2=\{3,4,5\},\ H_3=\{2,3\},\ H_4=\{1,4,6\}.$
	
	Feeding the transversal e-algorithm with $\HH_4$ yields (among others) the semifinal 01g-row $r$ in Table  3. It is good since it e.g. contains the minimal $\HH_4$-hitting set $\{1,2,5\}$. Yet $r$ is not very-good since $\{1,3,5\}\in r\setminus {\it MHS}(r)$ is a dud (viewing that $\{1,3\}\in {\it HS}(\HH_4)$).
	
	\begin{tabular}{l|c|c|c|c|c|c|c}
		& $1$ & $2$ & $3$ & $4$  &$5$  & $6$  &  \\ \hline
		& & & & & & &  \\ \hline
		$r=$ &  $g_1$ & $g_2$ & $g_2$ & $g_1$ & 1& $g_1$&  merely-good  \\ \hline
		& & & & & & &  \\ \hline
		$r_1=$ &  {\bf 1} & $g_2$ & $g_2$ & {\bf 0} & 1& {\bf 0}&  merely-good \\ \hline
		$r_2=$ &  {\bf 0} & $g_2$ & $g_2$ & {\bf 1} & 1& {\bf 0}&   very-good \\ \hline
		$r_3=$ &  {\bf 0} & $g_2$ & $g_2$ & {\bf 0} & 1& {\bf 1}&  merely-good \\ \hline\hline	
		$r_4=$ &  1 & 1 & 0 & 0 & 1& 0& very-good   \\ \hline
		$r_5=$ &  0 & 1 & 0 & 0 & 1& 1&  very-good \\ \hline	
		& & & & & & &  \\ \hline	
		$\rho_1=$ &  $g_1$ & {\bf 1} & {\bf 0} & $g_1$ & 1& $g_1$&  very-good  \\ \hline	
		$\rho_2=$ &  $g_1$ & {\bf 0} & {\bf 1} & $g_1$ & 1& $g_1$&  merely-good  \\ \hline\hline
		$\rho_3=$ & 0 &0 & 1 & 1 & 1& 0&  very-good  \\ \hline
	\end{tabular}
	
	{\sl Table 3: Replacing a merely-good row by very-good rows.}
	
	We strive to replace $r$ by disjoint rows which are very-good and jointly contain the same minimal hitting sets as $r$.
	It is natural to start by picking any g-wildcard of $r$, say $(g_1,g_1,g_1)$, and {\it expand} $r$ accordingly as $r=r_1\uplus r_2\uplus r_3$ (see Table 3). We call $r_1,r_2,r_3$ the {\it sons} of $r$.
	One checks that $\{2,4,5\},\{3,4,5\}\in {\it MHS}(\HH_4)$, and so $r_2$ is very-good. As to $r_1$, it is merely-good. Specifically, by expanding the second g-wildcard $(g_2,g_2)$ one obtains $r_1=(1,{\bf 0},{\bf 1},0,1,0)\uplus (1,{\bf 1},{\bf 0},0,1,0)$, where the first son  is bad since $\{1,3\}\in {\it HS}(\HH_4)$, and the second (call it $r_4$) is very-good. Also $r_3$ is merely-good; it decomposes as
	$r_3=(0,{\bf 0},{\bf 1},0,1,1)\uplus (0,{\bf 1},{\bf 0},0,1,1)$, where the first son is bad ($36\in {\it HS}(\HH_4)$) and the second (call it $r_5$) is very-good. To summarize, we managed to replace the semifinal merely-good row $r$ by the final very-good rows $r_2,r_4,r_5$.

	Alternatively, one can start by expanding $(g_2,g_2)$. This yields the rows $\rho_1,\rho_2$ in Table 3. One checks that $\rho_1$ is very-good, but $\rho_2$ is not. Specifically, when expanding $(g_1,g_1,g_1)$ in $\rho_2$, two of the three arising 01-rows are bad. The third one (labelled $\rho_3$) is very-good. To summarize, $r$ can even be replaced by {\it two} very-good rows, i.e. $\rho_1,\rho_3$.

 {\bf 6.2} The example above suggests the following method to replace a semifinal good row $r$  by  final very-good rows that jointly contain the same MHSes  as $r$. There is nothing to do if $r$ is already very-good. By induction assume that a stack is filled with disjoint merely-good 012g-rows which jointly contain exactly the  MHSes  contained in $r$. (Initially $r$ is the only member of the stack.) Remove the top row $r'$ from the stack. Expanding any g-wildcard of $r'$ yields candidate\footnote{There is no danger confusing the with the kind of candidate sons in 3.1.} sons $r_1,r_2,...$ akin to 6.1. The very-good  candidate sons are final; they are removed from the stack and stored somewhere else. The bad  ones are thrown away, and the merely-good  ones are put on top of the stack.It is clear that the new stack maintains the induction hypothesis. When the stack is empty, the final rows are disjoint
	and  jointly contain the same MHSes as $r$.

Above we used the wording "expanding any $g$-wildcard". Without going into details we mention that "any" needs not be random but can be chosen in ways that likely increase compression.

	\section{Deciding very-goodness using inclusion-exclusion}

 The larger our semifinal rows $\rho$ in (17) the more desirable is it to have efficient criteria for very-goodness and badness. In particular in Sec. 6  we reduced the handling of merely-good rows, to large extent,  to the existence of such tests. In this and the next Section we offer two very-goodness tests. The one in Section 7 relies on inclusion-exclusion.
	
	 {\bf 7.1} Consider a fixed semifinal 01g-row $\rho$ triggered by $\HH\s {\cal P}([w])$. We say that $Z\subseteq [w]$ is a {\it potential $\rho$-spoiler}  if there is a $Y\in \rho$ with $Z\uplus\{a\}=Y$. In Table 4  the set system of all potential $\rho$-spoilers of some semifinal $\rho$ is represented  as disjoint union $d_1\uplus\cdots\uplus d_5$ of 01g-rows. Its cardinality is $24+\cdots+6=74$. Generally the following holds:
	
\begin{itemize}	
\item[(20)] {\it With $\gamma_1,...,\gamma_t$ being the length of the $g$-wildcards of $\rho$, the number of potential $\rho$-spoilers of the semifinal row $\rho$ is $Pot=(\gamma_2\cdots\gamma_t)+(\gamma_1\gamma_3\cdots\gamma_t)+(\gamma_1\cdots
 \gamma_{t-1})+|ones(\rho)|\cdot\gamma_1\cdots\gamma_t$.}
 \end{itemize}

	\begin{tabular}{l|c|c|c|c|c|c|c|c|c|c|c|c|c|c|c|c|c|c|c|c}
		
		&1 &2 & &3 & &4 &5 &6 & & & &7 &8 &9 &10&  &11 &12& cardinality   \\ \hline
		
		& & &  & & & & & & & & & & & & & & & & & \\ \hline	
		$\rho=$ &  $g_1$ &  $g_1$ & & $0$& &  $g_2$ &  $g_2$ &  $g_2$& & &  &  $g_3$&  $g_3$ &  $g_3$ &  $g_3$ &  & $1$ &  $1$ & 630   \\ \hline
		
		& & & & & & & & & & & & & & & & & & & \\ \hline
		$d_1=$ &  $g_1$ &  $g_1$ & & $0$ & & $g_2$ &  $g_2$ &  $g_2$& & &  &  $g_3$&  $g_3$ &  $g_3$ &  $g_3$ & & ${\bf 0}$ &  $1$ & 24   \\ \hline
		
		$d_2=$ &  $g_1$ &  $g_1$ & & $0$ & & $g_2$ &  $g_2$ &  $g_2$& & &  &  $g_3$&  $g_3$ &  $g_3$ &  $g_3$ &  &$1$ &  ${\bf 0}$  & 24  \\ \hline
		
		$d_3=$ &  ${\bf 0}$ &  ${\bf 0}$ & & $0$ & & $g_2$ &  $g_2$ &  $g_2$& & &  &  $g_3$&  $g_3$ &  $g_3$ &  $g_3$ &  &$1$ &  $1$ & 12   \\ \hline
		
		$d_4=$ &  $g_1$ &  $g_1$ & & $0$ & & ${\bf 0}$ &  ${\bf 0}$ &  ${\bf 0}$ & & &  &  $g_3$&  $g_3$ &  $g_3$ &  $g_3$ &  &$1$ &  $1$  & 8  \\ \hline
		
		$d_5=$ &  $g_1$ &  $g_1$ & & $0$ & & $g_2$ &  $g_2$ &  $g_2$& & &  &  ${\bf 0}$&  ${\bf 0}$ &  ${\bf 0}$ &  ${\bf 0}$ &   &$1$ &  $1$ & 6   \\ \hline
		
		& & & & & & & & & & & & & & & & & & &\\ \hline
		$\delta_1=$ &  $g_1$ &  $g_1$ & & $0$ & & $1$ &  {\bf 0} &  {\bf 0}& &
		&  &  {\bf 0}&  {\bf 0} &  $g_3$ &  $g_3$ & & $0$ &  $1$ & 4   \\ \hline
		
		$\delta_2=$ &  $g_1$ &  $g_1$ & & $0$ & & $1$ &  {\bf 0} &  {\bf 0}& &
		&  &  {\bf 0}&  {\bf 0} &  $g_3$ &  $g_3$ & & $1$ &  $0$ & 4   \\ \hline
		
		$\delta_3=$ &  $0$ &  $0$ & & $0$ & & $1$ &  {\bf 0} &  {\bf 0}& &
		&  &  {\bf 0}&  {\bf 0} &  $g_3$ &  $g_3$ & & $1$ &  $1$ & 2   \\ \hline
		
		$\delta_4=$ &  $g_1$ &  $g_1$ & & $0$ & & $0$ &  {\bf 0} &  {\bf 0}& &
		&  &  {\bf 0}&  {\bf 0} &  $g_3$ &  $g_3$ & & $1$ &  $1$ & 4   \\ \hline
		
		$\delta_5=$ &  $g_1$ &  $g_1$ & & $0$ & & $1$ &  {\bf 0} &  {\bf 0}& &
		&  &  {\bf 0}&  {\bf 0} &  $0$ &  $0$ & & $1$ &  $1$ & 2   \\ \hline		
		
	\end{tabular}

	{\sl Table 4: Counting  $\rho$-spoilers by  applying inclusion-exclusion}
	
	For a  semifinal $\rho$ we define an {\it $\rho$-spoiler} as a potential $\rho$-spoiler that happens to be an $\HH$-hitting set. If $Sp=Sp(\rho,\HH)$ is the number of $\rho$-spoilers, then a moment's reflection confirms:
	
	(21)\quad The semifinal row $\rho$ is very-good  iff $Sp=0$.
	
	If say $H_i,H_j,H_\ell$ are hyperedges of $\HH$ then we define $N(i,j,\ell)$  as the number of potential $\rho$-spoilers $Z$ with $Z\cap H_i=Z\cap H_j=Z\cap H_\ell=\emptyset$.  Since a potential spoiler is a spoiler iff it cuts all hyperedges of $\HH$,
	we can compute  $Sp$  with inclusion-exclusion as
	
	(22)\quad $Sp=Pot-N(1)-N(2)-\cdots -N(h)+N(1,2)+\cdots+(-1)^h N(1,2,..,h)$.

	 Calculating $2^h$ terms $N(..)$ may seem inefficient but
	the larger $|\rho|$ and $w$, and the smaller $h$, the more inclusion-exclusion will prevail over
	the "naive" way in 5.2 which spends $O(hk^2)$ time per $k$-element member $X\in\rho$.
	
	{\bf 7.2}  Furthermore, based on the first three {\bf B}on{\bf f}erroni\footnote{These inequalities are the backrock of many theorems in probability theory and statistics. Consult any good textbook please.} inequalities  these implications often alleviate full-blown inclusion-exclusion:
	
	\begin{enumerate}
		\item[{(Bf1)}] $Pot-N(1)-\cdots -N(h)>0\ \Rightarrow\ Sp>0$ (not very-good)
		
		\item[{(Bf2)}] $Pot-N(1)-\cdots -N(h)+N(1,2)+\cdots +N(h-1,h)=0\ \Rightarrow\ Sp=0$  (very-good)
		
		\item[{(Bf3)}] $Pot-N(1)-\cdots +N(1,2)-N(1,2,3)-\cdots -N(h-2,h-1,h)>0$\\$ \Rightarrow\ Sp>0$  (not very-good)
	\end{enumerate}

{\bf 7.3} Full-blown inclusion-exclusion can also be avoided by other means. Recall that	
 $N(i_1,...,i_t)$ is the number of potential spoilers $Z$ with $Z\cap H_{i_1}=\cdots =Z\cap H_{i_t}=\es$. But this is equivalent to $Z\cap (H_{i_1}\cup\cdots \cup H_{i_t})=\es$. If the hyperedges are all very large (say of cardinality $>w/3$) then it is likely that $U:=H_{i_1}\cup\cdots\cup H_{i_t}=[w]$ even for small index sets $\{i_1,...,i_t\}\s [h]$. But then $N(i_1,...,i_t)=0$. (More generally "=0" happens iff $U$ contains a $g$-bubble or cuts $ones(\rho)$.)
 
 This appeals to the following more general endeavour (work in progress, arXiv:1309.6927v3). In every inclusion-exclusion problem the family of {\it relevant} index sets
  $\{i_1,...,i_t\}$, i.e. the ones  that satisfy $N(i_1,...,i_t)\neq 0$, constitute a set-ideal ${\cal N}\s {\cal P}([h])$. If this so-called {\it nerve} ${\cal N}$ is small and can be obtained in clever ways (i.e. not by scanning ${\cal P}([h])$ ), then inclusion-exclusion speeds up considerably.

	{\bf 7.4} According to (21) it follows from $Sp>0$ that $\rho$ is not very-good. But $\rho$ stays merely-good (as opposed to bad) unless $Sp$ sky-rockets. To make this more precise, let us generally order the sizes of the $g$-wildcards occuring in $\rho$ as $\gamma_1\le\cdots\le\gamma_t$. Then each $\rho$-spoiler $Z$ can prevent at most $\gamma_t$ many $X\in\rho$ from being in ${\it MHS}(\HH)$.  Since $|\rho|=(\gamma_1\cdots\gamma_{t-1})\gamma_t$, we conclude:
	
	(23)\quad If $\gamma_1\le\cdots\le\gamma_{t-1}\le\gamma_t$ and $Sp(\rho)<\gamma_1\gamma_2\cdots\gamma_{t-1}$, then $\rho$ is good.
	
	Although the bound $\epsilon_1\cdots\epsilon_{t-1}$ is sharp, in practise\footnote{Computational experiments have been carried out in a previous version arXiv:2008.08996v2 of this article. } it is likely that for much higher values of $SP(\rho)$
 the row $\rho$ remains merely-good.

\section{Deciding  very-goodness using Rado's Theorem}

Our second method to decide the very-goodness of a semifinal 01g-row $\rho_j$ is based on certain "critical" pairs $(\rho_i,\rho_j)$. Matroids [9] will play a crucial role. Let us jump into medias res with Rado's Theorem [9,p.702]:

\begin{itemize}
\item[(24)]  {\it Consider any matroid $M$ on a set $E$  and any family $\{Q_i:\ i\in I\}$ of subsets of $E$. Then  this family has a hitting set which is $M$-independent iff\\
\hspace*{0.8cm} $|J|\le rank\Big(\bigcup\{Q_j:\ j\in J\}\Big)$ for all $J\subseteq I$.}
\end{itemize}

{\bf 8.1} Apart from inviting matroids, here comes the second ingredient:

\begin{itemize}
\item[(25)]{\it  The semifinal row $\rho_j$ in (17) is not very-good iff there is a semifinal row $\rho_i\neq\rho_j$ such that $X\s Y$ for some $X\in\rho_i$ and $Y\in\rho_j$.}
\end{itemize}

{\it Proof of (25).} Assume that such $X,Y$ with $X\s Y$ exist. Since $X=Y$ is impossible ($\rho_i\cap\rho_j=\es$), we have $X\subset Y$. Since $Y$ properly contains a $\HH$-hitting set, we conclude $Y\not\in {\it MHS}(\HH)$, and so $\rho_j$ is not very-good. Conversely suppose that $\rho_j$ is not very-good. Picking any dud $Y\in\rho_j\setminus {\it MHS}(\HH)$ there is $X'\in {\it MHS}(\HH)$ with $X'\subset Y$. This $X'$ belongs to a unique semifinal row $\rho_i$ by (17). We have $\rho_i\neq\rho_j$ since $deg(\rho_i)=|X|<|Y|=deg(\rho_j)$. $\square$

In view of (25)   we call $(X,Y)$ a {\it spoiling pair for $\rho_j$} (not to be confused with the 'spoilers' in Sec. 7) if 

\begin{center}
$(Y\in\rho_j)\ \wedge\ (\exists i\neq j) (X\in\rho_i)\  \wedge\ X\subset Y.$
\end{center}

When $(X,Y)$ is a spoiling pair for $\rho_j$ then necessarily there is $i$ such that $(\rho_i,\rho_j)$ is a {\it critical} pair in the sense that $deg(\rho_)<deg(\rho_j)$ and $ones(\rho_i)\cap zeros(\rho_j)=\es$. This speeds up searching spoiling pairs $(X,Y)$ for likely-very-good rows $\rho_j$.

{\bf 8.2} To illustrate consider a hypothetical hypergraph that has triggered the two semininal rows $\rho_1,\rho_2$ in Table 5 . In fact $(\rho_2,\rho_1)$ is a critical pair since $deg(\rho_2)=4<5=deg(\rho_1)$ and $ones(\rho_2)\cap zeros(\rho_1)=\{5\}\cap \{1,2,3\}=\es$. In order to efficiently decide the existance of a spoiling pair $(X,Y)$ for $\rho_1$ (with $X\in\rho_i=\rho_2$), notice that any such $(X,Y)$ has $X\cap zeros(\rho_1)=\es$, and so $X\in \rho_2'$ (see Table 5). But why does $\rho_2'$ also differ from $\rho_2$ on the rightmost part? Because the $g_1g_1$ in $\rho_2$ was forced to become ${\bf 01}$.
Now ${\bf 1}$ in $\rho_2'$ triggers a $1$ at the same location in $\rho_1$, which  transforms $g_3g_3g_3$ in $\rho_1$ to $ 100$, i.e. replaces $\rho_1$ by $\rho_1'$.  Dropping the common 0's of $\rho_1',\rho_2'$ one gets two $1g$-rows $\rho_1'',\rho_2''$ with the same index set, in our case $E:=\{4,5,6,7,8,9,10,11,12\}$.

\begin{tabular}{l|c|c|c|c|c|c|c|c|c|c|c|c|c|c|c|c|c|c|c}
	&1 &2  &3  &4 & &5 & &6  &7 &8 & &9 &10 &11&  &12 &13 &14   \\ \hline
	
	& & & & & & & & & & & & & & & & & & \\ \hline

	$\rho_1=$ &  $0$ &  $0$ &  $0$ &  $1$ & & $1$ & & $g_1$ &$g_1$  & $g_1$& & $g_2$ & $g_2$ & $g_2$ & &$g_3$ & $g_3$ & $g_3$    \\ \hline
	
	$\rho_2=$ &  $g_1$ &  $g_2$ &  $g_3$ &  $g_2$ & & $1$ & & $g_2$ &$g_2$  & $g_3$ & & $g_3$ & $g_3$ & $g_2$ & &$g_1$ & $g_2$ & $g_2$    \\ \hline
	
	& & & & & & & & & & & & & & & & & & \\ \hline
	
	$\rho_1'=$ &  $0$ &  $0$ &  $0$ &  $1$ & & $1$ & & $g_1$ &$g_1$  & $g_1$ & &$g_2$ & $g_2$ & $g_2$ & &$1$ & $0$ & $0$    \\ \hline
	
	$\rho_2'=$ &  {\bf 0} &  $0$ &  $0$ &  $g_2$ & & $1$ & & $g_2$ &$g_2$  & $g_3$ & & $g_3$ & $g_3$ & $g_2$ & &{\bf 1} & $0$ & $0$    \\ \hline	
	
	& & & & & & & & & & & & & & & & & & \\ \hline
	
	$\rho_1''=$ &   &   &      &$1$ & & $1$&  & $g_1$ & $g_1$  & $g_1$ & &$g_2$ & $g_2$ & $g_2$ & &$1$ &  &     \\ \hline
	
	$\rho_2''=$ &   &   &      &$g_2$ & & $1$&  & $g_2$ & $g_2$   &$g_3$&  & $g_3$ & $g_3$ & $g_2$ & & $1$ &  &   \\ \hline	
		
\end{tabular}

{\sl Table 5: Deciding the existence of a spoiling pair with a Theorem of Rado}

 That's when the matroid takes over. Namely, the partition $E=\{4\}\uplus\{5\}\uplus\{6,7,8\}\uplus\{9,10,11\}\uplus\{12\}$ determined by the 1's and $g$-wildcards of $\rho_1''$ defines a so-called {\it partition matroid} $M=M(E)$ where by definition $X\subseteq E$ is $M$-independent iff $X$ cuts each part of the partition in at most one element. In contrast, the analogous partition induced by $\rho_2''$ is not used for a second matroid but rather yields the set system $\{Q_i:\ i\in I\}$ in (24). In our case $I=\{1,2,3,4\}$ and
$Q_1=\{4,6,7,11\},\ Q_2=\{5\},\ Q_3=\{8,9,10\},\ Q_4=\{12\}$. Consequently, if $X$ is an $M$-independent transversal of $\{Q_i:\ i\in I\}$, then $X$ extends to a spoiling pair $(X,Y)$ of $\rho_1''$ (and each spoiling pair arises this way). The existence of such spoiling pairs is handled by the rank condition in statement (24). Take say $J=\{2,3,4\}$. Then

$|J|=3\le 4=rank(Q_2\cup Q_3\cup Q_4)=rank(\{5,8,9,10,12\})$.

One sees that generally the cardinality of $I$ in (24) equals $deg(\rho_2'')$ which,  even for large hypergraphs $\HH$,  often is a modest number (and so all $J\subseteq I$ can be evaluated painlessly).

	\section{The benefits of having ${\it MC}(\HH)$ and   ${\it MinNotMC}(\HH)$ }

 In Section 9 we fill in gaps in 5.3 and 5.4 and deepen our understanding of the set-ideal ${\it MC}(\HH)$ of all MC-sets.  The acronym MC  [4] abbreviates { \it Minimality Condition}, i.e. the fact  that "being MC" is a necessary condition for "being a MHS". (So MC has nothing to do with {\bf MC}HS appearing in (13).) Subsection 9.1 proves the key fact  ${\it MHS}(\HH)={\it HS}(\HH)\cap {\it MC}(\HH)$. The set-system ${\it MinNotMC}(\HH)$ (consisting of the generators of the complementary set-filter of ${\it MC}(\HH)$) is introduced in 9.2. Our third criterion (after Sections 7 and 8) for very-goodness  appears in 9.3. In 9.4-9.5 we unravel the enigmatic set-systems $PotKi(\rho_i)$ from 5.3 and trim them to set-systems $Ki(\rho_i)\ (i\le R)$. Using Vertical Layout these $R$ set-systems can be calculated "simultaneously". 
 
 The remainder of Section 9 relies on the dual companion of the transversal $e$-algorithm, i.e. the noncover $n$-algorithm which we glimpsed in 3.3. In 9.6 the latter represents ${\it MC}(\HH)$ as a disjoint union (32) of 012n-rows. In a sense (32) dualizes (11).  The dualization continues in 9.7 in that 01g-rows get accompanied by $01g^*$-rows. Furthermore ${\it MHS}(\HH)$ is represented as disjoint union of set-systems
 $\rho_i\cap\sigma_j$, where the $\rho_i$'s are 01g-rows and the $\sigma_j$'s are $01g^*$-rows. Using inclusion-exclusion   $|\rho_i\cap\sigma_j|$ can be calculated quickly (9.8). This enables us to calculate $|{\it MHS}(\HH)|$ without knowing ${\it MHS}(\HH)$. Merely deciding whether or not $\rho_i\cap\sigma_j=\es$ works faster still and it e.g. leads to the badness-criterion (37).

{\bf 9.1} Given  $\HH\s {\cal P}([w])$, in 5.4 we defined MC-sets $X\s [w]$ and saw that the set-system ${\it MC}(\HH)$ of all MC-sets is a set-ideal. The intersection of this set-ideal with the set-filter ${\it HS}(\HH)$ turns out to be highly relevant:

 {\bf Theorem 2: }{\it For any   hypergraph $\HH$ it holds that  ${\it HS}(\HH)\cap {\it MC}(\HH)={\it MHS}(\HH).$}
	
 {\it Proof.} Take any $X\in {\it MHS}(\HH)$ and fix $b\in X$. There are $b$-hyperedges $H$, i.e. with $b\in H$, since otherwise $X\setminus\{b\}$ would remain a hitting set, in contradiction to $X$ being minimal. Suppose none of the $b$-hyperedges were to cut $b$ sharply from $X$. Then $(X\setminus\{b\})\cap H\neq\es$ for all $b$-hyperedges $H$, and of course $(X\setminus\{b\})\cap H'\neq\es$ for all other hyperedges $H'$. This contradicts $X\in {\it MHS}(\HH)$, and hence shows that $X\in {\it MC}(\HH)$. From ${\it MHS}(\HH)\s {\it HS}(\HH)$ follows ${\it MHS}(\HH)\s {\it HS}(\HH)\cap {\it MC}(\HH)$.

 Conversely pick $Y\in {\it HS}(\HH)\cap {\it MC}(\HH)$. Since by assumption $Y$ is a hitting set, it suffices to show that $Y\setminus\{b\}$ is no hitting set for all $b\in Y$. In view of $Y\in {\it MC}(\HH)$ some $H_0\in\HH$ cuts $b$ sharply from $Y$, hence  $H_0\cap ((Y\setminus\{b\})=\es$, hence $Y\setminus\{b\}$ is no hitting set. $\square$ 
 
 Consider any $X\in {\it HS}(\HH)\cap {\it MC}(\HH)$ and suppose $X$ was not a facet of ${\it MC}(\HH)$. Then there was a facet $Y$ with $X\subset Y$, and so $Y\in {\it HS}(\HH)\cap {\it MC}(\HH)$. But in view of Theorem 2 this yields  the contradiction of two comparable members of ${\it MHS}(\HH)$. We conclude that

 (26)\quad {\it At most the facets of ${\it MC}(\HH)$ can be minimal $\HH$-hitting sets.}

 In 5.4 we found that with respect to  $\HH_3$ the set $\{1,2\}$ is MC but no hitting set. One checks that $\{1,2\}$ is a facet of ${\it MC}(\HH_3)$. This shows that the by (26) necessary condition of being a  facet of ${\it MC}(\HH)$ is not sufficent for being an MHS.
	As another consequence of Theorem 2 we find that for each  semifinal 01g-row $\rho$ 
 from (17) we have
	
	(27) \quad $\rho\cap {\it MC}(\HH)=\rho\cap {\it HS}(\HH)\cap {\it MC}(\HH)=\rho\cap {\it MHS}(\HH).$ 
	
	Because ${\it MC}(\HH)\s {\cal P}([w])$ is a set-ideal by (18), we can consider the complementary set-filter ${\cal F}:={\cal P}([w])\setminus {\it MC}(\HH)$  (see 3.3).  This yields a neat {\it sufficient} condition for badness:

 (28)\quad  {\it If the semifinal 01g-row $\rho$ is such that $ones(\rho)$ is not MC, then $\rho$ is bad.}

 To prove it, all $X\in\rho$ are supersets of $ones(\rho)$, and so $ones(\rho)\in {\cal F}$ implies $X\in {\cal F}$. 
 
 A semifinal 01g-row satisfying (28) will be\footnote{If the 01g-row is a 01-row then of course {\it bad $\Leftrightarrow $ easy-bad}.} called {\it easy-bad}. A bad row which is not easy-bad is {\it sophisticated-bad}; an example was given in 5.3.1.
 
  {\bf 9.2} By definition the set-system \vspace{2mm}
  
 \centerline{${\it MinNotMC}(\HH)\hspace{1cm} (of\ cardinality\ mnMC)$} 
  
  consists of the generatos of the set-filter ${\cal F}$ in 9.1. To spell it out, ${\it MinNotMC}(\HH)$ consists of those subsets of $[w]$ which are not MC, but all their proper subsets are MC. While ${\it MinNotMC}(\HH)$ is beneficial, it is also expensive to compute. Before we turn to the benefits, here comes a toy example.
  
  {\bf 9.2.1}  It turns out (see Sec. 10) that
 ${\it MinNotMC}(\HH_2)$ is the set-system ${\cal G}$ in (8). Here $\HH_2$ is from (9). To summarize

 (29)\quad $\HH_2=\{125,34,456,135,26\}\ has\ {\it MinNotMC}(\HH_2)=\{123,15,126,256,134,345,346,246\}.$

 For instance, $Z:=\{2,5,6\}$ is {\it not} MC since no $\HH_2$-hyperedge cuts out $6$ sharply: $Z\cap\{4,5,6\}=\{5,6\}$ and $Z\cap\{2,6\}=\{2,6\}$. However, let us verify that all 2-subsets $Z'\subset Z$ (and hence all subsets)  {\it are} MC. For instance, take $Z'=\{5,6\}$. While still
 $Z'\cap\{4,5,6\}=\{5,6\}$, now $Z'\cap\{2,6\}$ works, i.e. equals $\{6\}$. Since also $Z'\cap\{1,2,5\}=\{5\}$, the set $Z'$ is MC. Similarly one checks that the other 2-subsets of $Z$, i.e. $\{2,6\}$ and $\{2,5\}$, are MC-sets.

  {\bf 9.3} We are now fit to return to the set-systems $PotKi(\rho)$ in 5.3. It follows from (27) that $X\in \rho$ is a dud iff $X$ is no MC-set, i.e. iff $X$ contains some member of ${\it MinNotMC}(\HH)$. In other words, setting 
  
  \centerline{$PotKi(\rho):={\it MinNotMC}(\HH)$}
  
  fulfils the requirement of 5.3 for whatever  semifinal 01g-row $\rho$. Trouble is, the set $PotKi(\rho)$ may be  bigger than it need be. Put another way, many members of $PotKi(\rho)$ are just  {\it potential} killers, i.e. dangerous for other rows, but not harming any  $X\in\rho$. Thus if $Z\in {\it MinNotMC}(\HH)$ is such that $Z\not\s X$ for all $X\in\rho$, we are led to say $Z$ is $\rho$-{\it harmless}.  Putting 
  
  $Harmless(\rho):=\{Z\in {\it MinNotMC}(\HH):\ Z\ is\ \rho {\it -harmless}\},$
  
  and
  
  $Ki(\rho):={\it MinNotMC}(\HH)\setminus Harmless(\rho)$

   we hence get a third very-goodness criterion:
   
  (30)\quad {\it A semifinal 01g-row $\rho$ is very-good iff $Ki(\rho)=\es.$}
  
  Recall that $\HH_2$ triggered the $R=4$ semifinal 01g-rows $\rho_1,...,\rho_4$ in (16), all of which happened to be of the same degree and hence very-good. In accordance with (30) one verifies that indeed $Ki(\rho_1)=\cdots=Ki(\rho_4)=\es$.
As crisp as (30) may look, viewing that ${\it MinNotMC}(\HH)$ is hard to find (Section 10), the criteria for very-goodness derived in Sections 7-8 remain attractive.

{\bf 9.4} The good news is, once ${\it MinNotMC}(\HH)$ has been conquered, VL will yield $Ki(\rho_i)$ {\it simultaneously} for all semifinal 01g-rows $\rho_i\ (1\le i\le R)$. Namely\footnote{What follows is only for VL-enthusiasts.}, we start by initializing certain auxiliary sets to $Ha(i):=\es$ for all $1\le i\le R$. For each fixed $Z\in {\it MinNotMC}(\HH)$ we will calculate the set $I(Z)$ of all  $i\le R$ which have $Z\in Harmless(\rho_i)$, and accordingly update $Ha(i):=Ha(i)\cup\{Z\}$ for all $i\in I(Z)$. Hence, once all $Z\in {\it MinNotMC}(\HH)$ have been processed, all $Ha(i)$ will have the correct content $Ha(i)=Harmless(\rho_i)$ (and so $Ki(\rho_i)={\it MinNotMC}(\HH)\setminus Ha(i)$ is obtained).

Calculating $I(Z)$ for fixed $Z$ works as follows. Say $\rho_1=(0,1,0,g_1,g_1,g_2,g_2,g_2,g_3,g_3)$. It will trigger the first three $01\infty$-{\it rows} $r_1,r_2,r_3$ of the matrix $A$ that underlies the VL application to come.  Turning all existing $1$'s to $0$'s,  setting all existing $0$'s to $\infty$ (more on that in moment),
and filling exactly one $g_i$-wildcard with $1$'s and the others with $0$'s, yields

$r_1=(\infty,0,\infty,1,1,0,0,0,0,0),\ r_2=(\infty,0,\infty,0,0,1,1,1,0,0),\ r_3=(\infty,0,\infty,0,0,0,0,0,1,1).$

In order to remember the {\bf n}umber $i=1$ of the {\bf s}emi{\bf f}inal 01g-row $\rho_i$  triggering $r_1,r_2,r_3$ we set $nsf(1)=nsf(2)=nsf(3):=1$. Say $\rho_2$ has two $g$-bubbles. Then it triggers analoguous
$01\infty$-rows $r_4,r_5$ (written below $r_1$ to $r_3$) and we record $nsf(4)=nsf(5):=2$. And so it goes on with $\rho_3$ up to $\rho_{R}$. 

Having calculated $A$ (say it has dimensions $41\times 10$), we can begin to process all $Z\in {\it MinNotMC}(\HH)$. If say $Z_1=\{5,6,7\}$, calculate the column $col:=col_5+col_6+col_7$. Then $col'=(1,2,0,...)$, where the fact that the second component is $\ge 2$ testifies that $Z_1$ cannot be contained in any member of the semifinal 01g-row with number $nsf(2)=1$ (since $Z_1$ cuts one $g$-bubble of that row in $\ge 2$ elements). As another example suppose that $Z_2$ is such that the corresponding length 41 column $col''$  has 20 components equal to 1, 20 equal to 0, and  the 13th component is $\infty$. How does this translate to plain language? It means that $Z_2$ is  harmless {\it only} for the semifinal 01g-row  $\rho_j\ (j:=nsf(13))$ because $Z_2\cap zeros(\rho_j)\neq\es$, and so $Z_2$ can't be contained in any member of $\rho_j$. (For all other semifinal 01g-rows $Z_2$ is a killer since it doesn't clash with their 0's and cuts all their g-bubbles in at most one element.)
For general $Z\in {\it MinNotMC}(\HH)$ with coupled column $col$ let $J$ is the position-set of the components $\ge 2$ that occur in $col$. By the above it is clear that $I(Z)=\{nsf(j):\ j\in J\}$ (it doesn't matter that $nsf(j)=nsf(j')$ for $j\neq j'$ is possible).

 {\bf 9.5}   In 5.3 two toy examples showed how $PotKi(\rho)$ helps to calculate $Duds(\rho)$. Let us propose a more systematic way 
 (from now on we drop $PotKi(\rho)$ and stick to $Ki(\rho)$) that in particular speeds up the detection of bad rows. So
put $Ki:=Ki(\rho)$. If $\rho$ has $t$ many g-bubbles then for all $1\le j\le t$ let $Ki[j]$ be the (possibly empty) set of $Z\in Ki$ that intersect exactly $j$ many g-bubbles (necessarily these intersections being singletons). Hence $Ki=Ki[0]\uplus Ki[1]\uplus\cdots\uplus Ki[t]$. Since each killer $Z\in Ki$ is necessarily disjoint from $zeros(\rho)$, we see that 
$Ki[0]=\{Z\in Ki:\ Z\s ones(\rho)\}$. Recalling the definition of "easy-bad" in 9.1 we claim:

(31)\quad $Ki[0]\neq\es$ {\it iff $\rho$ is easy-bad.}

{\it Proof of (31).} If $Z\in Ki[0]\neq\es$, then $Z$ (being a killer) is not-MC, hence the superset $ones(\rho)$ is not-MC, hence $\rho$ is easy-bad. If conversely $\rho$ is easy-bad, then $ones(\rho)$ (being not-MC) contains some $Z\in {\it MinNotMC}(\HH)$. Obviously $Z\in Ki[0]$. $\square$
     
   If $\rho$ is not easy-bad then it either is sophisticated-bad or 'actually-good'. To find out  fast, the  second most effective killers are the ones in $Ki[1]$  which we hence exploit first to inflate our changing set $Duds(\rho)$. Then we turn to $Ki[2]$, and so forth up to $Ki[t]$. If in the process $Duds(\rho)$ ever becomes $\rho$, we can stop and know that $\rho$ is bad. As seen in 5.3, VL can be used in all of that.
	

 {\bf 9.6} Recall from 3.3 that the noncover $n$-algorithm yields for each set-system ${\cal S}$ the family $NC({\cal S})$ of all ${\cal S}$-noncovers as a disjoint union of 012n-rows $\ol{\sigma_j}$. If in particular
 ${\cal S}:={\it MinNotMC}(\HH)$ then $NC({\cal S})={\it MC}(\HH)$. Therefore (11') specializes to 
 
 (32)\quad ${\it MC}(\HH)=\biguplus_{j=1}^{R^*}\ol{\sigma_j}$.

 For instance, recall that applying the transversal $e$-algorithm to $\HH_2$ yielded ${\it HS}(\HH_2)=\ol{\rho_1}\uplus\cdots\uplus\ol{\rho_4}$ (Table 2). If we dually apply the noncover $n$-algorithm to ${\it MinNotMC}(\HH_2)$ from (29) we get ${\it MC}(\HH_2)=\ol{\sigma_1}\uplus\cdots\uplus \ol{\sigma_6}$ (Table 6).

\begin{tabular}{l|c|c|c|c|c|c|c}
		& $1$ & $2$ & $3$ & $4$  &$5$  & $6$    & row-maximal sets \\ \hline
		& & & & & & &  \\ \hline
		$\ol{\sigma_1}=$ &  0 & 2 & $n$ & $n$ & $n$& 0&  $Max(\ol{\sigma_1})=\{234,235,245\}=(0,1,g^*,g^*,g^*,0)=:\sigma_1$  \\ \hline
		$\ol{\sigma_2}=$ &  0 & $n$ & 1 & 0 & $n$& 1&  $Max(\ol{\sigma_2})=\{236,356\}=(0,g^*,1, 0,g^*,1)=:\sigma_2$ \\ \hline
		$\ol{\sigma_3}=$ &  0 & 0 & 0 & 2 & 2& 1&  $Max(\ol{\sigma_3})=\{456\}=(0,0,0,1,1,1)=:\sigma_3$  \\ \hline
		$\ol{\sigma_4}=$ &  0 & 1 & 0 & 0 & 0& 1&  $Max(\ol{\sigma_4})=\{26\}=(0,1,0,0,0,1)=:\sigma_4$  \\ \hline
		$\ol{\sigma_5}=$ &  1 & 0 & $n$ & $n$ & 0& 2&  $Max(\ol{\sigma_5})=\{136,146\}=(1,0,g^*,g^*,0,1)=:\sigma_5$  \\ \hline
		$\ol{\sigma_6}=$ &  1 & 1 & 0 & 2 & 0& 0&  $Max(\ol{\sigma_6})=\{124\}=(1,1,0,1,0,0)=:\sigma_6$  \\ \hline
		
	\end{tabular}
	
	{\sl Table 6: Representing ${\it MC}(\HH_2)$ as disjoint union of  012n-rows}

 {\bf 9.7} Let us keep on dualizing. To begin with, for each $012n$-row $\ol{\sigma}$ in (32) one gets $Max(\ol{\sigma})$ by turning all 2's to 1's and all $n$-wildcards to $g^*$-wildcards, where by definition $(g^*,g^*,...,g^*)$ means "exactly one 0 here". For instance $\ol{\sigma_1}$ in Table 6 becomes $\sigma_1$. Generally each $\ol{\sigma_j}$ from (32) induces such a $01g^*$-{\it row} $\sigma_j$. Akin to (17) we claim that

 (33)\quad ${\it MHS}(\HH)\s\biguplus_{j=1}^{R^*}{\sigma_j}$.

 {\it Proof of (33).} From (32) and Theorem 2 follows ${\it MHS}(\HH)\s\biguplus_{j=1}^{R^*}\ol{\sigma_j}$. Hence each $X\in {\it MHS}(\HH)$
 is in a unique row $\ol{\sigma_j}$. We claim that $X\in Max(\ol{\sigma_j})=\sigma_j$. Indeed, since $X$ is a maximal member of ${\it MC}(\HH)$ by (26), it is a fortiori maximal within $\ol{\sigma_j}\s {\it MC}(\HH)$. $\square$

In  view of (33) we can carry over the concepts good, bad, very-good, and so on to $01g^*$-rows. For instance, as it is forced by (17) and (33), the 9 MHSes of $\HH_2$ appear both in (16) and Table 6, yet $R=4\neq 6=R^*$. Whereas all $\rho_i$ were very-good, $\sigma_4$ is bad; its only member $\{2,6\}$ is MC but no MHS. 

{\bf 9.7.1} Recall from Section 1 that our Main Quest is to retrieve the diamonds from the boxes (=semifinal 01g-rows) as efficiently as possible. As is evident from (33) one could also  retrieve the diamonds from dual boxes (= semifinal $01g^*$-rows). In fact, this is attempted in [4], yet in a one-by-one fashion based directly on ${\it MC}(\HH)$. We will continue to retrieve the MHSes from the semifinal 01g-rows but the "dual" $01g^*$-rows will play an important auxiliary role.
For technical reasons (see footnote in 9.8.1) the coupled 012e-rows and 012n-rows will resurface as well. For starters, observe that unless all involved rows are 01-rows it holds that $\rho_i\subset\ol{\rho_i}$ and  $\sigma_j\subset\ol{\sigma_j}$. Nevertheless, this takes place:
 
 \begin{itemize}
 \item[(34)] {\it For all $\ol{\rho_i},\rho_i$ in (11) and (17), and all $\ol{\sigma_j},\sigma_j$ in (32) and (33), we have $\ol{\rho_i}\cap \ol{\sigma_j}={\rho_i}\cap {\sigma_j}$.}
 \end{itemize}

{\it Proof of (34).} It suffices to show $\ol{\rho_i}\cap \ol{\sigma_j}\s{\rho_i}\cap {\sigma_j}$. As we long know, $\ol{\rho_i}\cap {\it MHS}(\HH)\s\rho_i$. Similarly, as shown above, $\ol{\sigma_j}\cap {\it MHS}(\HH)\s\sigma_j$. Together with Theorem 2 follows that
$\ol{\rho_i}\cap \ol{\sigma_j}=(\ol{\rho_i}\cap {\it HS}(\HH))\cap (\ol{\sigma_j}\cap {\it MC}(\HH))=(\ol{\rho_i}\cap {\it MHS}(\HH))\cap (\ol{\sigma_j}\cap {\it MHS}(\HH)\s {\rho_i}\cap {\sigma_j}$. $\square$

From Theorem 2, (11), (32), the distributivity of $\cap$ over $\cup$, and (34) follows 

$(35)\quad {\it MHS}(\HH)=(\biguplus_{i=1}^R \ol{\rho_i})\cap (\biguplus_{j=1}^{R^*} \ol{\sigma_j})
=\biguplus_{i,j\ge 1}(\ol{\rho_i}\cap \ol{\sigma_j}) = \biguplus_{i,j\ge 1}(\rho_i\cap \sigma_j).$

{\bf 9.7.2} To illustrate (35), taking $\ol{\rho_2}$ from Table 2 and $\ol{\sigma_1}$ from Table 6, it holds that $\ol{\rho_2}\cap\ol{\sigma_1}=\{235,245\} $. Generally speaking, intersecting 012e-rows with 012n-rows (or 01g-rows with $01g^*$-rows) is no easier than intersecting two 012e-rows (see 2.3). As one way out one can ponder to either expand the 012e-row or the 012n-rows as a disjoint union of 012-rows. For instance $\ol{\sigma_1}$ expands as shown in Table 7:

\begin{tabular}{l|c|c|c|c|c|c|}
		& $1$ & $2$ & $3$ & $4$  &$5$  & $6$     \\ \hline
		& & & & & &   \\ \hline
		$\ol{\sigma_1}=$ &  0 & 2 & $n$ & $n$ & $n$& 0  \\ \hline\hline
		$\ol{\sigma_{11}}=$ &  0 & 2 & ${\bf 0}$ & {\bf 2} & {\bf 2} & 0  \\ \hline
$\ol{\sigma_{12}}=$ &  0 & 2 & ${\bf 1}$ & {\bf 0} & {\bf 2} & 0  \\ \hline
$\ol{\sigma_{13}}=$ &  0 & 2 & ${\bf 1}$ & {\bf 1} & {\bf 0} & 0  \\ \hline
 	\end{tabular}
	
	{\sl Table 7: Expanding a 012n-row into 012-rows.}

It follows that $\ol{\rho_2}\cap\ol{\sigma_1}=(\ol{\rho_2}\cap\ol{\sigma_{11}})\uplus(\ol{\rho_2}\cap\ol{\sigma_{12}})\uplus(\ol{\rho_2}\cap
\ol{\sigma_{13}})$. Each term on the right, and generally each intersection of a 012e-row with a 012-row, is either empty or again a 012e-row (2.2.1). In our particular case $\ol{\rho_2}\cap\ol{\sigma_{11}}=(0,1,0,1,1,0),\ \ol{\rho_2}\cap\ol{\sigma_{12}}=(0,1,1,0,1,0),\ \ol{\rho_2}\cap\ol{\sigma_{13}}=\es$. 

 Let us argue that in the present scenario such intersections 
 are {\it always} either empty or 01-rows. So suppose $\ol{\rho}$ is from (11) and $\ol{\sigma}$
 from (32) got again expanded into 012-rows  $\ol{\sigma_*}$. Since each MHS $X$ contained in $\ol{\sigma}$ is maximal within $\ol{\sigma}$, it will also be maximal within the 012-row  $\ol{\sigma_*}$ it happens to lie. Any two MHSes being incomparable, there cannot be another MHS in $\ol{\sigma_*}$. Because $\ol{\rho}\cap\ol{\sigma_*}$
 ,if nonempty, is a 012e-row that by (35) consists entirely of MHSes, this 012e-row is actually a 01-row that matches $X$.

The bottom line is this. Formula (35) likely cannot  be exploited to compress ${\it MHS}(\HH)$; at most (35) can be used for one-by-one enumeration. Whether and when this competes with the four one-by-one schemes from Section 5, remains to be seen\footnote{In particular 5.4.1 will be stiff competition. As opposed to 9.2 to 9.8, the method in 5.4.1 does not even rely on ${\it MinNotMC}(\HH)$. }. For (35) to be competitive it will be necessary (possibly not sufficient) that all empty intersections $\rho_i\cap\sigma_j$ can be reckognized fast.

{\bf 9.8} The true calling of (35) is to find the {\it cardinality} $|{\it MHS}(\HH)|$! Namely, suppose  that: 

\begin{itemize}
\item[(36)] {\it For each $\ol{\rho_i}$ in (11) we can obtain (preferably few) 012n-rows
$\ol{\tau_1},...,\ol{\tau_m}$ such that\\ $\ol{\rho_i}\cap {\it MHS}(\HH)\s(\ol{\rho_i}\cap \ol{\tau_1})\uplus\cdots\uplus (\ol{\rho_i}\cap \ol{\tau_m})$ and $\ol{\tau_1},...,\ol{\tau_m}\s {\it MC}(\HH)$.}
\end{itemize}

In view of (35), statement (36) is plausible. A full proof of (36) that also touches on $Ki(\ol{\rho_i})$ and on computational issues will be given in Section 11.5.
Accepting (36) we first note that from $\ol{\rho_i}\cap \ol{\tau_j}\s {\it HS}(\HH)\cap {\it MC}(\HH)={\it MHS}(\HH)$ follows that "$\s$" in (36) in fact is "$=$". Hence\\ $|\ol{\rho_i}\cap {\it MHS}(\HH)|=|\ol{\rho_i}\cap \ol{\tau_1}|+\cdots +|\ol{\rho_i}\cap \ol{\tau_m}|$.
  Because $|{\it MHS}(\HH)|$ is the sum of  $R$ terms $|\ol{\rho_i}\cap {\it MHS}(\HH)|$, calculating $|{\it MHS}(\HH)|$ boils down to calculating
 $|\ol{\rho}\cap\ol{\sigma}|$ for an arbitrary 012e-row $\ol{\rho}$ and 012n-row $\ol{\sigma}$ (this problem occurs in other circumstances as well). Let us apply inclusion-exclusion to do so.
 
{\bf 9.8.1} To fix ideas take $\ol{\rho}:= (e_1, e_1, e_2, e_2, e_2, e_2)$
and $\ol{\sigma}:= (n_1, n_2, n_1, n_2, n_3, n_3)$. (The presence of entries 0,1,2 would only cause trivial changes
in the sequel.) 
Let $N(e_1),N(e_2),N(e_1e_2)$ be the numbers of bitstrings $x\in\ol{\sigma}$ that violate\footnote{Notice that $ee..e$ (of length $k$) is violated by just {\it one} bitstring $00..0$, whereas $gg..g$ is violated by $2^k-k$ many bitstrings. That's why in the context of inclusion-exclusion we prefer to deal with $|\ol{\rho}\cap\ol{\sigma}|$ rather than $|{\rho}\cap {\sigma}|$. Recall from (34) that $\ol{\rho}\cap\ol{\sigma}={\rho}\cap {\sigma}$.}, respectively, the $e_1$-bubble, the $e_2$-bubble, and both $e_i$-bubbles. By inclusion-exclusion it holds that

$|\ol{\rho}\cap\ol{\sigma}|=|\ol{\rho}|-N(e_1)-N(e_2)+N(e_1e_2)=27-12-4+1=12$

in view of $|\ol{\rho}|=3\cdot 3\cdot 3\ (see\ 2.2),\ N(e_1)=|(0,0,2,2,n_3,n_3)|=12,\ N(e_2)=|(2,2,0,0,0,0)|=4,\ N(e_1e_2)=|(0,0,0,0,0,0)|=1$.

Similarly (using obvious notation) one obtains 12 as 

$|\ol{\rho}\cap\ol{\sigma}|=|\ol{\sigma}|-N(n_1)-N(n_2)-N(n_3)+N(n_1n_2)+N(n_1n_3)+N(n_2n_3)-N(n_1n_2n_3)=45-16-16-12+4+4+4-1=12$

in view of $|\ol{\tau}|=3\cdot 15,...,N(n_3)=|(e_1,e_1,2,2,1,1)|=12,..,N(n_1n_2n_3)=|(1,1,1,1,1,1)|=1$.

In general we launch inclusion-exclusion on the row with the fewer wildcards.
  	Interestingly, deciding merely whether or not $\ol{\rho}\cap \ol{\tau}$ is {\it empty}, works even faster  than inclusion-exclusion (more on that in 11.5.2). 
 This speed of deciding the emptiness of $\ol{\rho}\cap \ol{\tau}$ prompts us to finally state a badness criterion  for semifinal 012e-rows. Thus, as opposed to (28), it is a sufficient {\it and necessary} (albeit somewhat clumsy) condition:

 \begin{itemize}
 \item[(37)] {\it Suppose the 012e-row $\ol{\rho_i}$ is as in (36).
  Then $\ol{\rho_i}$ is bad iff
 $\ol{\rho_i}\cap \ol{\tau_j}=\es$ for all $1\le j\le m$.}
 \end{itemize}

	\section{How to calculate ${\it MinNotMC}(\HH)$ in the first place}
	
	In order to understand how ${\it MinNotMC}(\HH_2)$ in (29) was computed\footnote{For various small $\HH\ (w<20)$ the author 
 ran both Algorithm 4 from 5.4 and the sophisticated method in spe to get ${\it MinNotMC}(\HH)$. Both always yielded the same result. While this is no formal proof of correctness of the method in spe, it is makes correctness very likely.}, it pays to momentarily relabel\footnote{Otherwise the elements of $H_1,...,H_6$ clash with the {\it indices} 1,2,...,6 of $H_1,...,H_6$.} the hyperedges of $\HH_2$ in obvious ways:
	
	(38)\quad $H_1=\{a,b,e\},\ H_2=\{c,d\},\ H_3=\{d,\epsilon,f\},\ H_4=\{a,c,\epsilon\},\  H_5=\{b,f\}$
	
	Let us refine the  property '$T$ is MC'. Thus for any set $T$ and fixed $u\in T$ we  say '$T$ {\it is $u$-MC}' if  $crit(u,T):=\{H\in\HH:\ H\cap T=\{u\}\}$ is nonempty. Consequently it holds for all $T\subseteq W:=\{a,b,c,d,\epsilon,f\}$ that:
	
	(39)\quad $T$ is MC \hspace*{0.3cm} $\Longleftrightarrow$ \hspace*{0.3cm} $T$ is $u$-MC for all $u\in T$
	
	(40)\quad $T$ is not-MC \hspace*{0.2cm} $\Longleftrightarrow$ \hspace*{0.2cm} $T$ is not-$u$-MC for some $u\in T\ \Longleftrightarrow\ crit(u,T)=\es$\hspace*{0.2cm} for some $u\in T$
	
	For instance $T=\{d,c,f\}$ is not-$d$-MC because from $d\in T\cap H_i$ always follows $|T\cap H_i|\ge 2$, the relevant indices being $i=2,3$. 
	
	{\bf 10.1} For $u\in W$ put ${\cal S}_u:=\{i\in [h]:\ u\in H_i\}$. Therefore

	(41)\quad ${\cal S}_a=\{1,4\},\ {\cal S}_b=\{1,5\},\ {\cal S}_c=\{2,4\},\ {\cal S}_d=\{2,3\},\ 
	{\cal S}_\epsilon=\{1,3,4\},\ {\cal S}_f=\{3,5\}.$ 
	
	The fact hat $\{d,c,f\}$ is not-$d$-MC can now be seen as tantamount to ${\cal S}_d \subseteq
	{\cal S}_c \cup {\cal S}_f$. Generally the not-$u$-MC sets bijectively match the {\it set coverings} of ${\cal S}_u$
	by {\it other}  ${\cal S}_v$'s.
	
	Our aim is to calculate the family ${\it MinNotMC}(\HH_2)$ of minimal not-MC sets. According to (40) they are found {\it among} the minimal not-$u$-MC sets, where $u$ ranges over $W$. Let us hence find for each fixed $u\in W$ all minimal set coverings of ${\cal S}_u$. The systematic method follows in 10.2, but for $\HH=\HH_2$ we can proceed by inspecting (41):
	
	\begin{itemize}
		\item The minimal set coverings of ${\cal S}_a$ are $\{{\cal S}_b,{\cal S}_c\},\ \{{\cal S}_\epsilon\}$
		
		\item The minimal set coverings of ${\cal S}_b$ are 
		$\{{\cal S}_a,{\cal S}_f\},\ \{{\cal S}_\epsilon,{\cal S}_f\}$
		
		\item The minimal set coverings of ${\cal S}_c$ are 
		$\{{\cal S}_a,{\cal S}_d\},\ \{{\cal S}_d,{\cal S}_\epsilon\}$
		
		\item The minimal set coverings of ${\cal S}_d$ are 
		$\{{\cal S}_c,{\cal S}_\epsilon\},\ \{{\cal S}_c,{\cal S}_f\}$
		
		\item The minimal set coverings of ${\cal S}_\epsilon$ are 
		$\{{\cal S}_a,{\cal S}_d\},\ \{{\cal S}_a,{\cal S}_f\},\
		\{{\cal S}_b,{\cal S}_c,{\cal S}_d\},\ \{{\cal S}_b,{\cal S}_c,{\cal S}_f\}$
		
		\item The minimal set coverings of ${\cal S}_f$ are 
		$\{{\cal S}_b,{\cal S}_d\},\ \{{\cal S}_b,{\cal S}_\epsilon\}$	     
	\end{itemize}
	
	Therefore the minimal not-$a$-MC sets are $\{{\bf a},b,c\}$ and $\{{\bf a},\epsilon\}$, and so forth until the 
	minimal not-$f$-MC sets are $\{{\bf f},b,d\}$ and $\{{\bf f},b,\epsilon\}$. The inclusion-minimal sets among these sets\footnote{ For instance $\{\epsilon,a,d\}$ gets killed by $\{a,\epsilon\}$, and e.g. the double occurence of $\{b,\epsilon,f\}=\{f,b,\epsilon\}$ can be pruned. In general proceed as in 2.5 to get the minimal sets.} are (in shorthand notation)
	$abc,a\epsilon,baf,b\epsilon f,cad,cd\epsilon,dcf,fbd$.
	Relabelling back $a\to 1,...,f\to 6$ yields ${\it MinNotMC}(\HH_2)$ in (29).
	
	{\bf 10.2} As is well known, finding minimal set coverings is cryptomorphic to finding minimal hypergraph transversals. Let us make this cryptomorphism explicite by recalculating the  set coverings of the set ${\cal S}_\epsilon$ by the set system ${\cal S}:=\{{\cal S}_a,{\cal S}_b,{\cal S}_c,{\cal S}_d,{\cal S}_f\}$. Because ${\cal S}_\epsilon=\{1,3,4\}$, at least one member of ${\cal S}$ must cover 1; only ${\cal S}_a,{\cal S}_b$ can do that. Similarly, only ${\cal S}_d,{\cal S}_f$ can contain 3, and only ${\cal S}_a,{\cal S}_c $ can contain 4. Thus we define the {\it auxiliary hypergraph} triggered by $\epsilon$ as
	
	$(42)\ \HH_2^{aux}(\epsilon):=\{  \{{\cal S}_a,{\cal S}_b \},\{{\cal S}_d,{\cal S}_f \},\{{\cal S}_a,{\cal S}_c \} \}.$
	
	It follows that the $\HH_2^{aux}(\epsilon)$-hitting sets are exactly the minimal set-coverings of ${\cal S}_\epsilon$ by other ${\cal S}_u$'s.
	It is natural to employ again the transversal e-algorithm to calculate all minimal $\HH_2^{aux}(\epsilon)$-hitting sets.
	
	The transversal e-algorithm starts by imposing the hyperedge $\{{\cal S}_a,{\cal S}_b \}$ of $\HH_2^{aux}(\epsilon)$, and then imposes $\{{\cal S}_d,{\cal S}_f \}$. Since the two happen to be disjoint, this is achieved by the single 012e-row $r_1$ in Table 7. Imposing $\{{\cal S}_a,{\cal S}_c \}$ upon $r_1$ yields the two final rows $r_2$ and $r_3$.
	It happens that both of them are very-good, i.e. $Min(r_2)$ and $Min(r_3)$ need not  be pruned further.
	
	\begin{tabular}{l|c|c|c|c|c|c}
		& ${\cal S}_a$ & ${\cal S}_b$ & ${\cal S}_c$ & ${\cal S}_d$  & ${\cal S}_f$& $Min(r_i)$ \\ \hline
		& & & & & &   \\ \hline
		$r_1=$ &  $e$ & $e$ & $2$ & $e'$ & $e'$&     \\ \hline
		& & & & & &  \\ \hline
		$r_2=$ &  {\bf 1} &  2 &{\bf 2} & $e'$ & $e'$&  
		$\{\{{\cal S}_a,{\cal S}_d\},\ \{{\cal S}_a,{\cal S}_f\}\}$  \\ \hline
		$r_3=$ &   {\bf 0} & 1  &{\bf 1} & $e'$ & $e'$&   
		$\{\{{\cal S}_b,{\cal S}_c,{\cal S}_d\},\ \{{\cal S}_b,{\cal S}_c,{\cal S}_f\}\}$ \\ \hline	
	\end{tabular}
	
	{\sl Table 7: Calculating all the minimal set coverings of ${\cal S}_\epsilon$ with the $e$-algorithm}

	\section{Numerical experiments}

While terminology and overall structure of the article in front of you have improved a lot compared to the 2021-version (arXiv:2008.08996v2), there is a problem concerning the 2021 Mathematica experiments: The author lacked time to substitute some of the 2021 subroutines by implementations of the superior ideas discussed in previous Sections. After long deliberation I decided to pick a few of the most telling numerical experiments done in 2021, recast them in Table 8 below, and describe them thoroughly with adapted terminology (i.e. from the present article). All experiments are characterized by the {\it signature} $(w,h,k)$ that refers to a hypergraph $\HH\s {\cal P}[w]$ whose $h$ hyperedges $H\in\HH$ are random and have {\it uniform cardinality} $|H|=k$. For some signatures (in 11.1) we managed to calculate $|{\it MHS}(\HH)|$ exactly. For other signatures $|{\it MHS}(\HH)|$ could only be {\it approximated}; e.g. because ${\it MinNotMC}(\HH)$ could not be conquered (11.2), or not even the $R$ semifinal rows in (11) could (11.3). In 11.4 and 11.5 we speculate on  future improvements. Finally 11.6 compares our "wildcard-approach" with an algorithm of Toda [5] which is based on BDD's and which therefore also offers some kind of compression..

{\bf 11.1} Whenever $|{\it MHS}(\HH)|$ could be determined exactly, the procedure usually was as follows. The (transversal) $e$-algorithm, fed with $\HH$, terminates and outputs $R$ many semifinal 01g-rows $\rho_i$ (see (17)). Whenever ${\it MinNotMC}(\HH)$ could be calculated, then likewise  all $R$ set-systems $Ki(\rho_i)$ could be calculated (though not yet with the nifty VL way of 9.4). In this situation the potential very-goodness of $\rho_i$ (and if yes,  $|\rho_i|$) is settled at once in view of (30).  How to process the remaining merely-good or bad rows $\rho_j$? We mostly used 9.5 (inflating $Duds(\rho_j)$ by processing $Ki[0]\uplus\cdots\uplus Ki[t]$) or\footnote{The author does not remember for each signature occuring in the 2021 experiments which variant was used. Notice that the 9.6 - 9.8 variant is more powerful but has an overhead, i.e. 9.5 might be faster for moderate $|\rho_j|$. }
 9.6 - 9.8 (combining the $n$-algorithm with inclusion-exclusion)

Thus one hypergraph $\HH$ of signature (60,20,5) (see Table 8) triggered $R=26701$ semifinal 01g-rows $\rho_i$  of average degree 13. The calculation took 13 seconds.
	 Calculating ${\it MinNotMC}(\HH)$ of cardinality $mnMC=309$ took 1 second. Using the 9.6-9.8  way we found that $\HH$ had 51'109'682 MHSes. Perhaps more informative than knowing $|{\it MHS}(\HH)|$ is it to know the average 1914 of the (absolute) contents $|\rho_i\cap {\it MHS}(\HH)|$, as well as the average {\it relative} content $|\rho_i\cap {\it MHS}(\HH)| / |\rho_i|=0.77\ (77\%)$. (Up to small rounding error one retrieves $|{\it MHS}(\HH)|$ by multiplying with $R$ the average absolute  content.) As to the (30,50,70)-hypergraph, since its semifinal rows have little content and ${\it MinNotMC}(\HH)$ is large, the 9.5 way was faster. 
 For some  (70,20,30)-hypergraph the precise value of $|{\it MHS}(\HH)|$ was obtained {\it without} the aid ${\it MinNotMC}(\HH)$ because Algorithm 4 from 5.4 managed to process all semifinal rows  (including the very-good-ones) one-by-one.

{\bf 11.2} For some $(w,h,k)$-hypergraphs $\HH$ it was possible to calculate all $R$ semifinal rows but not the exact value of ${\it MHS}(\HH)$. That is because either ${\it MinNotMC}(\HH)$ was too hard to calculate (see also 11.4) and Algorithm 4  not up to the task. Or, while ${\it MinNotMC}(\HH)$ {\it could} be obtained, either  $R$ or the sizes $|\rho_i|$ were too large to process, in whatever way, the not very-good rows (see also 11.5). In this situation we picked 1000 among the $R$ semifinal rows at random\footnote{In fact we picked the first thousand 012e-rows produced by the transversal $e$-algorithm. One may object that these rows are not  representative because they match the
’leftmost’ 1000 nodes of the computation tree. This objection can be dismissed as follows. Whenever a top row of the LIFO stack (3.1) gets removed, we switch the new top
row with a random row further below. The effect is that the first 1000 semifinal rows are as random as any other sample of 1000 semifinal rows.   } and used them to approximate the average content of semifinal rows. 

There is one $\HH$ which doesn't quite fit "In this situation".
For this $\HH$ of signature (100,40,3) the 113 potential killers in ${\it MinNotMC}(\HH)$ could be calculated in just 0.4 sec. Among the 10367 semifinal rows 94\% were very-good (identified via $Ki(\rho_i)=\es$) and their  cardinalities summed up to 3190986028403520327. The remaining semifinal rows were all merely-good and still very dense. The 9.5 variant being out of question due to the size of $\rho_i$, the author  speculates (but doesn't remember fore sure) that attempting the 9.6 - 9.8 variant must have failed due to the inferior 2021 subroutine for inclusion-exclusion (see 11.5.2). 

{\bf 11.3} In some cases not  all semifinal rows could be generated, i.e. the $e$-algorithm failed and $R$ was unknown. Nevertheless, one can still employ the $e$-algorithm to generate {\bf 1000} random semifinal 01g-rows. The last three lines in Table 8 arose this way. It is interesting to compare the signatures (100,40,3) and (100,80,3), as well as  (30,50,7) to (30,5000,7). As usual, if $w,k$ stay fixed while $h$ increases, the absolute content "deteriorates".

As to the last column in Table 8, if all $R$ semifinal 012e-rows (and whence semifinal 01g-rows) could be classified (whether or not $|{\it MHS}(\HH)|$ was achieved) then we evidently get the exact percentages of very-good, merely-good, and bad rows. They appear (rounded) in the last column. If not all semifinal rows could be computed, then the numbers in the last column were extrapolated by applying the Monte Carlo method  to the 1000 semifinal rows (be it in 11.2 or 11.3) that {\it were} computed .

 {\bf 11.4}   As to calculating ${\it MinNotMC}(\HH)$, considerably less time was spent for running the $w$ many auxiliary
	transversal e-algorithms than for minimizing the resulting set system ${\cal S}$ to $Min({\cal S})={\it MinNotMC}(\HH)$. For instance for the (30,50,7)-instance it took only 61 seconds to calculate ${\cal S}$ (of cardinality 252'211), but 2503 seconds to shrink
	${\cal S}$ to ${\it MinNotMC}(\HH)$ (of cardinality 55538). For the (70,20,30)-instance ${\it MinNotMC}(\HH)$ could not be calculated in reasonable time. Problem is, the minimization method used was inferior to the  ideas in 2.5 and 2.6.1.  
 
 {\bf 11.4.1} Is there hope compressing ${\it MHS}(\HH)$ without knowing ${\it MinNotMC}(\HH)$? Yes there is: While the inclusion-exclusion method of Section 7 has been experimented with in the 2020-version of the present article (arXiv:2008.08996v1, Section 6.3) only for small values $h$, there is hope (recall  7.3) to trim it considerably. Also Rado's Theorem (Section 8) should be kept in mind as basis for a very-goodness criterion.
 When trying to compress ${\it MHS}(\HH)$, Section 6 plays an important role as well. A second look at Section 6 shows that one can handle matters with just a very-goodness criterion. Put another way, a nontrivial badness-criterion is nice-to-have but not strictly necessary to repackage merely-good rows into fresh very-good rows.
 
 {\bf 11.5} As to merely calculating $|{\it MHS}(\HH)|$, let us first prove (36) from 9.8. There are two approaches to obtain the required 012n-rows $\ol{\tau_i}$. Both are based on possessing ${\it MinNotMC}(\HH)$. The first approach obtains the rows $\ol{\sigma_j}\ (j\le R^*)$ in (32) by feeding the whole of ${\it MinNotMC}(\HH)$ to the (noncover) $n$-algorithm. Then for all $1\le j\le R^*$ we check whether or not $\ol{\rho_i}\cap\ol{\sigma_j}=\es$ (see 11.5.2) and take as 
 $\{\ol{\tau_1},...,\ol{\tau_{m_i}}\}$ the set of all $\ol{\sigma_j}$ with $\ol{\rho_i}\cap\ol{\sigma_j}\neq\es$. The second approach only feeds $Ki(\ol{\rho_i})$ instead of 
 ${\it MinNotMC}(\HH)$ to the $n$-algorithm and thus obtains $m_i'$ many 012n-rows $\tau_j'$ that also do the job.

 {\bf 11.5.1} What are the pros and cons of the two approaches to provide each semifinal 012e-row $\ol{\rho_i}$ with "its" 012n-rows guaranteed by (36)? For starters, while the calculation of the $R$ set-systems $Ki(\ol{\rho_i})$  is {\it based} on ${\it MinNotMC}(\HH)$, it works smoothly with Vertical Layout (9.4). Since all $Ki(\ol{\rho_i})$ are small subsets of ${\it MinNotMC}(\HH)$, applying
 the $n$-algorithm to a single $Ki(\ol{\rho_i})$ takes much less time than applying it to ${\it MinNotMC}(\HH)$. Under circumstances even the {\it sum} of all $Ki(\ol{\rho_i})$-times compares well to the ${\it MinNotMC}(\HH)$-time; e.g. when $R$ is small and/or many $\ol{\rho_i}$ are very-good due to $Ki(\ol{\rho_i})=\es$ and hence  need not undergo the $n$-algorithm. How does $m_i$ compare to $m_i'$? In lockstep with the shorter time also the number $m_i'$ of produced 012n-rows will be smaller than the corresponding number $m_i$. Finally observe that by construction all $m_i$ rows $\ol{\tau_j}$ intersect $\ol{\rho_i}$, whereas this need not be the case for the $m_i'$ many rows $\tau_j'$. If, while running the $n$-algorithm on $Ki(\ol{\rho_i})$, one keeps on discarding candidate sons $\tau'$ with $\ol{\rho_i}\cap\tau'=\es$ (see 11.5.2), then it is guaranteed that no final 012n-row will be disjoint from $\ol{\rho_i}$. In this way one can further reduce $m_i'$ but perhaps that's not worth the effort. More computational experiments need to be carried out to clarify all of 11.5.1.

 {\bf 11.5.2} Two more loose ends must be addressed. First, the type of inclusion-exclusion proposed in 9.8.1 for calculating $|\ol{\rho}\cap \ol{\tau}|$ is superior to the type of inclusion-exclusion employed in the 2021-experiments of Table 8. Namely, as detailed in [arXiv:2008.08996v1, Sec.7.2], this slower kind of inclusion-exclusion relies on a bipartite graph whose shores are the e-wildcards of $\ol{\rho}$ and the n-wildcards of $ \ol{\tau}$  respectively. Since in 9.8.1 we only need one kind of wildcards, the 9.8.1 implementation in spe is up to $2^t$ times faster that the current implementation (where $t$ is the number of the kind of wildcards of which there are more).
 
 Second, deciding merely whether or not $\ol{\rho}\cap \ol{\tau}$ is {\it empty} works  faster still than 9.8.1 type  inclusion-exclusion. For starters, the intersection is clearly empty when 1's in one row clash with 0's in the other row. However, there can be more hidden reasons for emptiness; e.g. $(1,1,e,e)\cap (n_1,n_2,n_1,n_2)=\es$. The gory details of deciding the emptiness of $\ol{\rho}\cap \ol{\tau}$ have been tackled in [arXiv:2008.08996v1, Sec.8], yet all of that will be recast in a separate publication that also relates the matter to deciding the satisfiability of certain Boolean functions (of type $Horn\wedge AntiHorn$). Another issue is the Mathematica implementation of it all, and its possible overhead that slows it down for small size inputs.

	
	\begin{tabular}{|c|c|c|c|c|}
		$(w,h,k)$	&$R$, av.deg  & mnMC & content (abs/rel) & vg, mg, bad       \\ \hline
		& & & &  \\ \hline
$(60,20,5)$ & 26701, 13 (13s) & 309, (1s) & $1914,\ 77\%\ (5042s)$ & 43,49,8    \\ \hline
$(70,20,30)$ & $77448,5\ (39s)$ & --- & $8.5,\ 20\%\ (186s)$ & 13,62,25   \\ \hline
$(30,50,7)$ & $123584,10\ (56s)$ & 55538 (2564s) & $1.02,\ 20\%\ (66s)$ & 15,26,59  \\ \hline\hline

 $(70,20,5)$ & $13577,14\ (8s)$ & 256, (3s)& $113116,\ 86\%\ $ & 68,32,0,   \\ \hline
 $(70,20,6)$ & $41319,12\ (21s)$ & 730, (3s) & $1694,\ 82\%\ $ & 37,62,2,   \\ \hline
$(70,20,12)$ & $917377,10\ (1546s)$ & --- & $42,\ 27\%$ & 	33,60,7  \\ \hline
 $(100,40,3)$ & $10367,33\ (13s)$ & 113, (0.4s) & $3\cdot 10^{10},\ 99\%$ & 94,6,0 \\ \hline\hline

$(30,5000,7)$ & ${\bf 1000},18\ (103s)$ & --- & $0.45,\ 22\%$ & $23,14,64$   \\ \hline
$(100,80,3)$ & ${\bf 1000},34\ (2s)$ & 437, (2s) & $3000,\ 29\%$ & $6,73,21$   \\ \hline
$(10000,100,1000)$ & ${\bf 1000},14\ (158s)$ & ---  & $10^6,\ 77\%$ & $55,40,5$\\ \hline
		
	\end{tabular}
	
	{\sl Table 8.  Numerical evaluation and extrapolation of the minhit algorithm}

	{\bf 11.6} In [1] nineteen methods to calculate ${\it MHS}(\HH)$ have been pitted against each other on a common platform, using a variety of real-life datasets. 
	Our method\footnote{ As we have seen, we proposed many variants to achieve various subgoals, but for brevity let us stick with "our method". (In a previous version of this article we bothered to name all these variants.)} does not post factum fit that platform. For one thing, it is implemented in high-level Mathematica code and so far only ran on the author's laptop (Dell Latitude 7410). Furthermore, much different from [1], all hypergraphs in Table 8 have random and equicardinal hyperedges (which in view of 11.3 may be disadvantageous). Nevertheless, let us attempt a preliminary comparison  with two specific algorithms investigated in [1]. First, the {\it Murakami-Uno-algorithm [4]} (like us to some extent)) relies  on the MC-condition but proceeds one-by-one. Second, building on ideas of Knuth, the {\it [Toda-algorithm [5]}, like us, uses compression, but in  more implicite ways (BDD's). These two algorithms also happen to be the champions\footnote{Let us cite from [1, Sec. 5.4]:  {\it The algorithms of [4] and [5] are far faster than their competitors across a variety of input set families. Toda is extremely fast on inputs for which it terminates (...). However it frequently exhausted the 32GB available memory on our workstation.}
		In [1,Sec. 5.4] the authors describe a collection of 128 million minimal transversals as '{\it enormous}'. This is true, but only if listed one-by-one. } in [1]. Since the MC-condition has received  plenty attention in Sections 9 and 10, let us devote the remainder of 11.6 to the [5]-algorithm. Here come four aspects where our method seems to win out (but since talk is cheap  only direct confrontation can ultimately determine the pros and cons of both). 
	
	\begin{itemize}
		\item [(i)] As is well known (and repeated in [6]), having the BDD of a Boolean function $f$ yields at once the cardinality of the model set $Mod(f)$. With a bit more effort (but in linear total time) one gets the model set of $f$ as a disjoint union of 012-rows. Unfortunately, when $Mod(f)={\it MHS}(\HH)$, then the models are mutually incomparable, and so all 012-rows are necessarily 01-rows, i.e. no compression is achieved. (This is akin to 9.7.2.) Matters are  alleviated but not cured by Toda's use of zero-supressed BDD's (=ZDD's). Thus the ZDD provides an implicite compression of ${\it MHS}(\HH)$ which often provided $|{\it MHS}(\HH)|$ faster than the 18 competitors in [1]. But since ${\it MHS}(\HH)$ is only\footnote{If all hitting sets are encoded in a BDD then it follows from arXiv:1703.08511 that, if not the minimal, the transversals of fixed cardinality can be output in a compressed format (using $g$-wildcards).}   output one-by-one this didn't always mean overall victory.  The Toda-algorithm is probably faster than us whenever the compression-rate\footnote{The compression rate $|{\it MHS}(\HH)|/R$ not only depends on the structure of $\HH$ but also on the {\it order} in which the hyperedges are imposed. No research in that direction has been undertaken.} is low, such as for the (30,5000,7) signature. With increasing compression-rate the tables begin to turn. Also keep in mind: Our more pleasantly compressed  representation of ${\it MHS}(\HH)$ may be desirable enough that spending extra time on it is worthwile.
		
		\item [(ii)] Even when the final BDD is moderate in size, intermediate BDD's can be excessively large, thus causing memory problems. In contrast, the LIFO stack used by the transversal $e$-algorithm can never contain more than $h$ rows (this is  a classic result about LIFO stacks).
		
		\item [(iii)] In [5,p.101] Toda hopes to eventually parallelize one part of his algorithm, i.e. the calculation of a BDD that captures ${\it HS}(\HH)$. In contrast, parallelizing\footnote{Why parallelization (aka distributed computing) works smoothly in all LIFO-scenarios is e.g. explained in [6,Sec. 6.5]. }  our equivalent (the transversal e-algorithm)  is straightforward. In fact, the evaluation of all semifinal 012e-rows can be parallelized as well.
		
		\item [(iv)] Like our method some algorithms in [1] have the potential for cut-off (4.3), but the Toda-algorithm seems not to be  among them since it does not appear in Table 9 or 10 of [1].
	\end{itemize}

	\section{Enumerating all exact hitting sets }

	In  our last Section all our hypergraphs  $\HH\subseteq {\cal P}([w])$ of cardinality $h:=|\HH|$  are {\it full} in the sense that $\bigcup\HH=[w]$ (to avoid trivial cases). 
	An {\it exact hitting set (EHS)} with respect to a hypergraph $\HH$ is a subset $X\subseteq [w]$ such that $|X\cap H|=1$ for all $H\in\HH$.
	Because of $\bigcup\HH=[w]$  each $a\in X$ belongs to some hyperedge $H$. This implies that each EHS $X$ is\footnote{In the sense that for {\it each} $a\in X$ {\it every} $H$ containing $a$ cuts it out sharply.}  "very MC", and so a minimal hitting set. The converse fails\footnote{In fact $\HH_2$ in 3.1 has {\it no} EHSes. We mention in passing that hypergraphs $\HH'$ with ${\it MHS}(\HH')=E{\it HS}(\HH')$ can be reckognized in polynomial time [10], and that ${\it MHS}(\HH')$ can be output one-by-one with polynomial delay. The most obvious instance of ${\it MHS}(\HH')=E{\it HS}(\HH')$ occurs when the hyperedges of $\HH'$ are mutually disjoint.}. 
	
	In the sequel we compress the set $E{\it HS}(\HH)$ of all exact $\HH$-hitting sets  by 'imposing' the hyperedges one after the other (12.2-12.3). In doing so the previously used 01g-cards will be applicable even more directly, yet the trivial feasibility test (10) becomes much harder. One consequence (12.4) concerns the enumeration of all perfect matchings in certain graphs. Sections 12.1 and 12.5 deal with a  natural (apparently novel) equivalence relation induced on $[w]$ by every hypergraph $\HH\s{\cal P}([w])$. It prompts one to distinguish 'degenerate' and 'nondegenerate' hypergraphs.

 {\bf 12.1}  For a hypergraph $\HH=\{K_1,\ldots,K_h\}\subseteq {\cal P}(W)$ we say that $x,y\in W$ are ($\HH$-){\it equivalent}\\ (written $\sim$) if $\forall (1\le i\le h)\ x\in K_i\Leftrightarrow y\in K_i$. If the equivalence relation $\sim$ is the identity relation, then $\HH$ is called {\it nondegenerate}, otherwise {\it degenerate}. For instance, if $\HH$ is the  hypergraph of all stars of a graph (see 12.4) then $\HH$ is  nondegenerate. On the other hand, the vertices $6,8$ are $\HH_1$-equivalent (see 2.4.1), and so $\HH_1$ is degenerate.
	For each index set
	$I\subseteq [h]$ let $\HH(I)$ be the set of $a\in W$ which are in all $K_i$'s $(i\in I)$ and nowhere else. Formally
	
	(43)\quad $\HH(I):=\bigcap \{K_i:\ i\in I\}\ \cap\ \bigcap\{W\setminus K_i:\ i\in [h]\setminus I\}$.
	
	If $\HH(I)\neq\emptyset$ then $\HH(I)$ is a $\sim$class, and each  $\sim$class arises this way\footnote{ Once more VL can be used. In brief, letting $A$ be the $h\times w$ whose $i$th row is the characteristic bitstring of the $i$th hyperedge, one checks that $1\sim k$ iff ($BitOr(col[1],col[k])=col[1]$ and $BitAnd(col[1],col[k])=col[1]$). In this way the $\sim$class $\ol{1}$ can be determined. Next pick any $j\in [w]\setminus \ol{1}$ and determine $\ol{j}$ likewise. And so forth. }. It follows that $2^h<w$ is a sufficient condition for $\HH$ to be degenerate.

 \begin{itemize}
 \item[(44)] {\it Let $\HH$ be a hypergraph and let $r$ be any 01g-row contained in $E{\it HS}(\HH)$. Then each $g$-bubble $\{a,b,...\}$ of $r$ is contained in a $\sim$class. }
 \end{itemize}

{\it Proof of (44).} Let $K\in\HH$ be arbitrary with $a\in K$. By symmetry it suffices to show that $b\in K$. By way of contradiction suppose $b\not\in K$. Fix any $X\in r$ with $a\in X$ (by definition of $01g$-row there is such $X$). Then $X\cap K=\{a\}$ since $X$ is an {\it exact} hitting set. If $Y$ arises from $X$ by switching $a$ with $b$ then still $Y\in r$. But $Y\cap K=\es$, which contradicts the fact that $Y$ (being in $r$) is an (exact) hitting set. $\square$

	{\bf 12.2} Consider the hypergraph $\HH_5\s {\cal P}[9]$ consisting of the three hyperedges
	
	(45) \quad $K_1 =\{2,3,4,6\},\ K_2=\{1,2,3,4,5,7\},\ K_3=\{2,8,9\}$.
	
	If instead of $\{K_1,K_2,K_3\}$ we just have $\{K_1\}$, then the set of $\{K_1\}$-hitting sets, i.e.\\
	$\{X\s [9]: |X\cap \{2,3,4,6\}|=1\}$, can be written\footnote{ In Section 12 our familiar 01g-rows must be slightly generalized to 012g-rows.} as the 012g-row $r_0$ below.
	
	\begin{tabular}{l|c|c|c|c|c|c|c|c|c|c}
		
		&1 &2 &3 &4 &5 &6 & 7&8 &9&    \\ \hline
		
		& & & & & & & & &  &  \\ \hline
		
		$r_0=$ &  2 &  $g$ &  $g$ &  $g$ &  2 &  $g$ &   2&  2 &  2 & pending $K_2$     \\ \hline
		& & & & & & & & &  &  \\ \hline	
		$r'_1=$ &  2 &  ${\bf 0}$ &  ${\bf 0}$ &  ${\bf 0}$ &  2 &  1 &    2&  2 &  2 &  \\ \hline
		$r'_2=$ &  2 &  ${\bf g}$ &  ${\bf g}$ &  ${\bf g}$ &  2 &  0 &    2&  2 &  2 &  \\ \hline\hline
		
		$r_1=$ &  ${\bf g}$ &   0 &   0 &   0 &  ${\bf g}$&  1 &    ${\bf g}$&  2 &  2 & pending $K_3$ \\ \hline
		$r_2=$ &  {\bf 0} &  ${ g}$ &  ${ g}$ &  ${ g}$ &  {\bf 0} &  0 &    {\bf 0}&  2 &  2 & pending $K_3$ \\ \hline\hline
		& & & & & & & & &  &  \\ \hline
		
		$r_3=$ &  $g_1$ &   0 &   0 &   0 &  $g_1$&  1 &    $g_1$&  $g_2$ &  $g_2$ & final  \\ \hline
		$r_2=$ &  0 &  ${ g}$ &  ${ g}$ &  ${ g}$ &  0 &  0 &    0&  2 &  2 & pending $K_3$ \\ \hline
		& & & & & & & & &  &  \\ \hline
		
		$r_4=$ &  0 &  0 &  ${ g_1}$ &  ${ g_1}$ &  0 &  0 &    0&  $g_2$ &  $g_2$ & final \\ \hline
		$r_5=$ &  0 &  1 &  0 &  0 &  0 &  0 &    0&  0 &  0 & final \\ \hline			
	\end{tabular}
	
	{\sl Table 9: The working stack for the g-algorithm}
	
	In order to sieve the $\{K_1,K_2\}$-EHSes $X$ from $r_0$  we observe that $K_1\cap K_2=\{2,3,4\}$ and accordingly write $r_0=\gamma_1\uplus\gamma_2$ (Table 9). That helps because sieving the $\{K_1,K_2\}$-EHSes from the auxiliary rows $r_1',r_2'$  is easy. It results in $r_1,r_2$ respectively. For both rows the imposition of $K_3$ is still {\it pending}. Each row in the  stack must be tagged with this kind of information. Picking the top row of the current {\it working stack} $\{r_1,r_2\}$ we focus on $r_1$. It is evident that the subset of all $X\in r_1$ with $|X\cap K_3|=1$ can be written as the 012g-row $r_3$ in Table 9.  Row $r_3$ is {\it final} in the sense that all hyperedges have been imposed on it; this amounts to $r_3\s E{\it HS}(\HH)$. We hence remove $r_3$  from the  working stack and make it the first final row. It is clear that imposing $K_3$ on the last row $r_2$ in the working stack yields the final rows $r_4,r_5$. We hence have $EHT(\HH_1)=r_3\uplus r_4\uplus r_5$. In particular $\HH_5$ has 
	$3\cdot 2+2\cdot 2+1=11$ exact hitting sets.

	{\bf 12.3} In order to generally impose a hyperedge $K$ upon a 012g-row we erect a certain Abraham-flag (boldface in Table 10) akin to (7). Thus imposing $K=\{1,2...,6\}$ upon the 012g-row\footnote{To avoid distraction we often choose $twos(r)=ones(r)=zeros(r)=\emptyset$. Only trivial modifications would occur otherwise.} $\ol{r_0}$ in Table 10 yields $\ol{r_1}$ to $\ol{r_4}$.
	
	\begin{tabular}{l|c|c||c|c||c||c|c|c|c|c|c|c|c|}
		
		&1 &2 &3 &4 &5 &6 & &7 &8 &9 &10 &11 &12  \\ \hline
		
		& & & & & & & & & & & & &   \\ \hline
		
		$\ol{r_0}=$ &  $g_1$ &  $g_1$ &  $g_2$ &  $g_2$ &  $g_3$ &  $g_4$ &  &  $g_1$&  $g_1$ &  $g_2$ &  $g_3$ &  $g_3$ &  $g_4$    \\ \hline
		
		& & & & & & & & & & & & &   \\ \hline
		
		$\ol{r_1}=$ &  ${\bf g_1}$ &  ${\bf g_1}$ &  ${\bf 0}$ &  ${\bf 0}$ &  ${\bf 0}$ &  ${\bf 0}$ &  &  $0$&  $0$ &  $1$ &  $g_3$ &  $g_3$ &  $1$    \\ \hline
		
		$\ol{r_2}=$ &  ${\bf 0}$ &  ${\bf 0}$ &  ${\bf g_2}$ &  ${\bf g_2}$ & ${\bf 0}$  & ${\bf 0}$ &  &  $g_1$&  $g_1$ &  $0$ &  $g_3$ &  $g_3$ &  $1$    \\ \hline
		
		$\ol{r_3}=$ &  ${\bf 0}$ &  ${\bf 0}$ & ${\bf 0}$ & ${\bf 0}$ &  ${\bf 1}$ & ${\bf 0}$ &  &  $g_1$&  $g_1$ &  $1$ &  $0$ &  $0$ &  $1$    \\ \hline
		
		$\ol{r_4}=$ &  ${\bf 0}$ &  ${\bf 0}$ &  ${\bf 0}$ & ${\bf 0}$ &  ${\bf 0}$ &  ${\bf 1}$ &  &  $g_1$&  $g_1$ &  $1$ &  $g_3$ &  $g_3$ &  $0$    \\ \hline		
	\end{tabular}
	
	{\sl Table 10: Imposing the exact hitting set $\{1,\ldots,6\}$ upon the row $\ol{r_0}$}
	
	 Adhering to the terminology of 3.1 we call $\ol{r_1}$ to $\ol{r_4}$ the {\it candidate sons} of $\ol{r_0}$ (that arise upon imposing $K$ on $\ol{r_0}$). Again we need to know which of the candidate sons $\ol{r_i}$ are {\it feasible} in the sense that $\ol{r_i}\cap E{\it HS}(\HH)\neq \emptyset$, and infeasible candidate sons (=duds)
	should be cancelled.
  The popular Dancing-Links algorithm of Knuth which decides (though not in polynomial time) whether or not a given hypergraph admits a hitting set, is easily adapted to a feasibility test for candidate sons. Again the surviving candidate sons of $\ol{r_0}$ are called its {\it sons}. The described method will be coined\footnote{This name was previously used by the author in other circumstances involving g-wildcards. There is no danger of confusion.} the {\it g-algorithm}.

	{\bf Theorem 3:} {\it Let $\HH\subseteq{\cal P} [w]$ be a hypergraph. Then $E{\it HS}(\HH)$ can be enumerated as a disjoint union of $R$ many 01g-rows in time $O(Rhw\cdot feas(h,w))$. Here $feas(h,w)$ upper-bounds the time for any chosen subroutine (e.g. Dancing-Links) to decide whether a hypergraph with $\le w$ vertices and $\le h$ hyperedges has an EHS.}
	
	{\it Proof.} Throughout the g-algorithm the  top rows in the LIFO-stack  match the nodes of a computation tree (rooted at $r_0$) whose $R$ leaves are the final rows. The length of a root-to-leaf path equals the number of impositions that were required to generate that leaf (=final row), and hence that length is at most $h$. In the worst case (i.e. when all root-to-leaf paths are mutually disjoint and have maximal length) the number of non-root nodes, i.e. the number of impositions, equals $Rh$. 
	
	What is the maximum cost $imp(h,w)$ of imposing a hyperedge on a LIFO top row $r$?
	Building the at most $\tau=\tau(\HH):=max \{ |H|:\ H\in\HH\}$ 
	candidate sons $r_i$ of $r$ (by way of  0g0-Abraham-Flags) costs $O(\tau w)$.   Letting $feas(h,w)$ be any time bound\footnote{For technical reasons we postulate that $f(h,w)\ge hw$. For every non-trivial hypergraph $\HH$ this inequality holds anyway.} for checking  the feasibility of a 012g-row it costs $O(\tau feas(h,w))$ to discard the infeasible candidate sons. A surviving son $r_i$ satisfies a fixed hyperedge $K$ iff in $r_i$ the bits with indices in $K$ are all  0's except for one 1. Hence it costs $O(\tau hw)$ to tag each son with its pending hyperedge.
	We conclude that $imp(h,w)=O(\tau w+\tau hw+\tau feas(h,w))$ and
	therefore:
	
	(46)\quad The overall cost of imposing the hyperedges of $\HH$ in order to pack all exact hitting sets of\\ \hspace*{1cm}
	$\HH$ into $R$ disjoint 01g-rows  is $O(Rh\cdot imp(h,w))=O(Rh\tau(hw+feas(h,w)) )$.
	
	Since we postulated $f(h,w)\ge hw$ and since $\tau\le w$, we have $O(Rh\tau(hw+feas(h,w)) )=O(Rhw\cdot feas(h,w)).\ \square$

	{\bf 12.4} An important  kind of exact hitting set arises from any  graph $G$ with vertex set $V$ and edge set $E$. Namely, if $star(v)$ is the set of all edges incident with vertex $v$ and $\HH:=\{star(v):\ v\in V\}\subseteq {\cal P}(E)$, then the EHSes of $\HH$ are exactly the perfect matchings of $G$. Recall that $K_{3,3}$ is the complete bipartite graph both shores of which having 3 vertices. The bipartite graph with 3 vertices on each shore, such that each vertex is adjacent to every vertex on the opposite shore, is commonly denoted as $K_{3,3}$. A graph $G$ is $K_{3,3}$ {\it minor-free} if one cannot obtain $K_{3,3}$ from $G$ by deleting edges and vertices of $G$, nor by contracting edges of $G$. 
	
	{\bf Theorem 4:} {\it All perfect matchings of a $K_{3,3}$ minor-free graph $G$ can be enumerated in polynomial total time.}
	
	{\it Proof.} In our context each feasibility test performed by the $g$-algorithm  on a $01g$-row $r$ de facto decides whether a certain
	minor $G(0,1)$ of $G$ of has a perfect matching. Specifically, the 0's in  $r$ delete edges from $G$ which thus becomes a sparser graph $G(0)$. The 1's in $r$ constitute a partial matching $P$ in $G(0)$ which wants to be extended to a perfect matching of $G(0)$. This is possible iff a certain subgraph $G(0,1)$ of $G(0)$ has a perfect matching. Namely, $G(0,1)$ is obtained by removing all edges  of $P$, along with all edges incident with them. The arising isolated vertices are also removed.
	With $G$ also its minor $G(0,1)$ is $K_{3,3}$-minor-free. By Corollary 1 in [11] one can decide in polynomial time (in fact even NC-time) whether $G(0,1)$ has a perfect matching. Hence the function $feas(h,w)$ in Theorem 1 is bound by a polynomial in $h,w$, causing the overall algorithm to run in total polynomial time. $\square$
	
	One can dispense with $K_{3,3}$-minor-freeness if one allows for randomization because deciding the existence of a perfect matching is in RNC [12,p.347]. Perfect matchings in bipartite graphs have been dealt with before [13].

{\bf 12.5}	Let $\HH=\{K_1,...,K_h\}\s {\cal P}([w])$ be a hypergraph. Generally, if a $\sim$class $C$ intersects $K_i$, then it must be {\it contained} in $K_i$; otherwise there were $x,y\in C$, one in $K_i$, the other not, which is impossible. Therefore, if
	$\ol{K_i}$ denotes the set of $\sim$classes contained in $K_i$, then $K_i=\biguplus \ol{K_i}$. 
	The {\it reduced} hypergraph $\ol{\HH}:=\{\ol{K_1},\ldots,\ol{K_h}\}$ has $h_0\le h$ hyperedges and is nondegenerate. For instance, for $\HH_5$ in (45) the $\HH_5$-classes are  
	$\ol{1}(=\ol{5}=\ol{7})=\{1,5,7\},\ \ol{2}=\{2\},\ \ol{3}=\{3,4\},\ \ol{6}=\{6\},\ \ol{8}=\{8,9\}$. Hence 
	$\ol{\HH_5}=\{\ol{K_1},\ol{K_2},\ol{K_3}\}$, where $\ol{K_1}=\{\ol{2},\ol{3},\ol{6}\},\ \ol{K_2}=\{\ol{1},\ol{2},\ol{3}\},\ \ol{K_3}=\{\ol{2},\ol{8}\}$.
	
	 Let us connect all of this with $g$-wildcards. The $g$-bubble of the $g$-wildcard in row $r_0$ of Table 9 is $\{2,3,4,6\}$.
	Since this  is just $K_1$, it is a union of $\sim$classes. It follows at once from induction and the design Abraham-Flags that this property gets perpetuated:
	
	\begin{itemize}
		\item[(47)] {\it When applying the g-algorithm to the hypergraph $\HH$,  each occuring g-bubble is a union of $\sim$-classes.}
	\end{itemize}
	
	However, once the $g$-algorithm has terminated, all final 01g-rows are subsets of $E{\it HS}(\HH)$, and so by (44) all their $g$-bubbles are contained in single $\sim$classes.
This is compatible with (47) only if
  each $g$-bubble of a final row actually {\it is} an $\HH$-class.
	
	{\bf 12.5.1} In particular, when applying the $g$-algorithm to a nondegenerate hypergraph, each final $01g$-row must be a $01$-row (=bitstring).
	For instance, applying the  $g$-algorithm to the nondegenerate hypergraph $\ol{\HH_5}$ would  give the final $01$-rows in the left part of Table 11:
	
	\begin{tabular}{l|c|c|c|c|c|c |c|c|c|c|c|c |c|c|c|c}
		& $\ol{1}$ & $\ol{2}$ & $\ol{3}$ & $\ol{6}$  & $\ol{8}$&    & 1&5&7&2&3&4&6&8&9&\\ \hline
		& & & & & &  & & & & & & & & & &  \\ \hline
		&  1 & 0 & 0 & 1& 1&$\Ra$  &$g_1$ &$g_1$ &$g_1$ &0 &0 &0 &1 &$g_2$ &$g_2$ &    \\ \hline
		
		&  0 & 0 & 1 & 0& 1&$\Ra$  &0 &0 &0 &0 &$g$ &$g$ &0 &$g$ &$g$ &   \\ \hline
		&  0 & 1 & 0 & 0& 0&$\Ra$  &0 &0 &0 &1 &0 &0 &0 &0&0 &   \\ \hline
	\end{tabular}
	
	{\sl Table 11: The $g$-algorithm necessarily enumerates $E{\it HS}(\ol{\HH_5})$ one-by-one}
	
	One retrieves the final 01g-rows on the right in Table 11 by inflating each $1$ at position $\ol{k}$ on the left to a $g$-wildcard as large as the
	class $\ol{k}$ (with the understanding that $1$ stays $1$ if $\ol{k}$ is a singleton).
	
	{\bf 12.6} What is the bottom line in all of that? A devil's advocate might argue: For nondegenerate hypergraphs the  $g$-algorithm offers no compression, and for degenerate hypergraphs $\HH$ the compression can also be achieved by enumerating the EHSes of $\ol{\HH}$ with any {\it other algorithm}, and then inserting $g$-wildcards in a trivial manner.
	
	Here is the defender's argument: As elementary as they are, the  concepts  'degenerate' and 'nondegenerate' are new.
	Likewise for '$g$-wildcards' and 'Abraham-Flags'. Concerning 'other algorithm',  the author could not google any publication concerning the enumeration of all exact hitting sets of a  general hypergraph. Even concerning specific hypergraphs, the algorithm in [13] seems to be the  only publication.  
 
	{\bf 12.6.1} What is the importance of "degenerate/or not" in the context of ${\it MHS}(\HH)$? As testified by $\HH_2$ in (9), the MHSes of {\it nondegenerate} hypergraphs are often compressible nevertheless. For degenerate $\HH$ one could, as we did for EHSes, run all our techniques on the reduced hypergraph $\ol{\HH}$ and later compress further. Whether that actually gives better compression than just sticking to $\HH$ remains to be seen.

 {\bf 12.7 Conclusion:} This article promotes the compression of ${\it MHS}(\HH)$ by the use of wildcards. This approach is very promising for sparse hypergraphs (see (100,40,3) in Table 8), but not advisable for dense ones (see (30,50,7) in Table 8). As observed already in [3], what works particularly well in  the sparse case (which we henceforth assume) is the compression of $MC{\it HS}(\HH)$, i.e. of the minimum-cardinality hitting sets. As to compressing the remainder ${\it MHS}(\HH)\setminus MC{\it HS}(\HH)$, we apologize for having overwhelmed (or not?) the reader with a plethora of topics: Three criteria for very-goodness, many uses of Vertical Layout, the fact that ${\it MHS}(\HH)={\it HS}(\HH)\cap {\it MC}(\HH)$, the proposal and calculation of ${\it MinNotMC}(\HH)$, the primal-dual approach ($e$- and $n$-wildcards) for finding $|{\it MHS}(\HH)|$, and more. While often illustrated with luscious toy-examples, many of these ideas  await  implementation and comparison with other approaches (collaboration is welcome). The author also appreciates to be informed of further (some are given in [6]) real-life examples of hypergraphs with few but large hyperedges. As a "side show" Section 12 turned to exact (as opposed to minimal) hitting sets. The issue of when $E{\it HS}(\HH)$ can be compressed is more clear-cut (12.6) than it was for ${\it MHS}(\HH)$. Further we touched upon Knuth's Dancing-Links and on compressing all perfect matchings of a graph.

	\section*{References}
	\begin{enumerate}

		\item[{[1]}]  A. Gainer-Dewar, P.  Vera-Licona, The minimal hitting set generation problem: algorithms and computation. SIAM J. Discrete Math. 31 (2017), no. 1, 63-100.

  \item[{[2]}] T. Eiter, G. Gottlob, K. Makino, New results on monotone dualization and generating hypergraph transversals. SIAM J. Comput. 32 (2003), no. 2, 514-537.

  \item[{[3]}]  M. Wild, Counting or producing all fixed cardinality transversals. Algorithmica 69 (2014), no. 1, 117-129.
		
		\item[{[4]}] K. Murakami, T. Uno, Efficient algorithms for dualizing large-scale hypergraphs. Discrete Appl. Math. 170 (2014), 83-94.

		\item[{[5]}] T. Toda, Hypergraph transversal computation with binary decision diagrams, in: SEA 2013 Rome, Italy. 

  \item[{[6]}] M. Wild, ALLSAT compressed with wildcards: From CNF's to orthogonal DNF's by imposing
		the clauses one by one, The Computer Journal, Vol.65 (2022) 1073-1087.

  \item[{[7]}] M. Wild, J. Svante, S. Wagner, D.Laurie, Coupon collecting and transversal of hypergraphs, Discrete Mathematics an Theoretical Computer Science 2013, 259-270.

   \item[{[8]}] L. Shi, X. Cai, An exact fast algorithm for minimum hitting set, 2010 Third International Conference on Comp Sc. and Optimization.

    \item[{[9]}] A. Schrijver, Combinatorial Optimization, Algorithms and Combinatorics 24, Springer-Verlag Berlin Heidelberg 2003.

    	\item[{[10]}] T. Eiter,
		Exact transversal hypergraphs and application to Boolean $\mu$-functions. (English summary)
		J. Symbolic Comput. 17 (1994), no. 3, 215-225.	

  	\item[{[11]}]  V. Vazirani, NC algorithms for computing the number of perfect matchings in K3,3-free graphs and related problems. Inform. and Comput. 80 (1989), no. 2, 152-164.

  	\item[{[12]}] R. MotwaniP. Raghavan, Randomized Algorithms, Cambridge University Press 1995.

		\item[{[13]}] T. Uno, A fast algorithm for enumerating bipartite perfect matchings. Algorithms and computation (Christchurch, 2001), 367-379, Lecture Notes in Comput. Sci., 2223, Springer, Berlin, 2001.

			

	\end{enumerate}

\end{document}